\documentclass{article}
\usepackage{amsmath}
\numberwithin{equation}{section}
\usepackage{amsfonts,amssymb,latexsym,amsbsy, bbm, theorem,enumerate,color}
\usepackage{graphicx, graphics}
\usepackage{float,color,fancybox,shapepar,setspace,hyperref}
\usepackage{mathrsfs,amscd,amssymb,bm,psfrag,subfigure,harpoon,esvect}
\usepackage{subfigure}
\usepackage{pgf,tikz}
\usetikzlibrary{arrows}
\newtheorem{Main Theorem}{Main Theorem}

\newtheorem{Theorem}{Theorem}

\newtheorem{Lemma}{Lemma}

\newtheorem{Observation}{Observation}
\newtheorem{Claim}{Claim}

\newtheorem{definition}{Definition}
\def\square{\hbox{\vrule height8pt depth0pt
\vbox{\hrule width7.2pt\vskip7.2pt\hrule width7.2pt}\vrule
height8pt depth0pt}\smallskip}

\def\pf{\medskip\noindent {\emph{\bf Proof}.}~~}

\newcommand\ex{\ensuremath{\mathrm{ex}}}

\newcommand\cF{{\mathcal F}}

\newcommand\cN{{\mathcal N}}

\newcommand\cT{{\mathcal T}}

\usepackage{authblk}
\begin{document}

\title{The maximum number of triangles in $F_k$-free graphs}

\author{Xiutao Zhu$^{1,2,}$\footnote{zhuxt@smail.nju.edu.cn} , Yaojun Chen$^{1,}$\footnote{ yaojunc@nju.edu.cn} , D\'aniel Gerbner$^{2,}$\footnote{gerbner@renyi.hu}, Ervin Gy\H{o}ri$^{2,}$\footnote{Corresponding  author. Email:  gyori.ervin@renyi.hu},
Hilal Hama Karim$^{3,}$\footnote{hilal.hamakarim@edu.bme.hu}\\
 \small{$^1$Department of Mathematics, Nanjing University, Nanjing 210093, P.R. CHINA}\\
 \small{$^2$Alfr\'ed R\'enyi Institute of Mathematics, Budapest}\\
 \small{$^3$Department of Computer Science and Information Theory, Faculty of Electrical Engineering and Informatics, Budapest University of Technology and Economics, Budapest}
 }
\date{}
\maketitle
\begin{abstract}
\vskip 2mm
The generalized Tur\'{a}n number $ex(n,K_s,H)$ is the maximum number of complete graph $K_s$ in an $H$-free graph on $n$ vertices. Let $F_k$ be the friendship graph consisting of $k$ triangles. Erd\H{o}s and S\'os (1976) determined the value of $ex(n,K_3,F_2)$.
 Alon and Shikhelman (2016) proved that $ex(n,K_3, F_k)\le (9k-15)(k+1)n.$
In this paper, by using a method developed by Chung and Frankl in hypergraph theory, we determine the exact value of $ex(n,K_3,F_k)$ and the  extremal graph for any $F_k$ when $n\ge 4k^3$.
\vskip 5mm
\noindent{\bf Keywords}: Generalized Tur\'{a}n number, triangle, friendship graph

\end{abstract}

\section{Introduction}

One of the most basic problems in extremal Combinatorics is the study of the \textit{Tur\'an number} $\ex(n,F)$, that is the largest number of edges an $n$-vertex $F$-free graph can have. A natural generalization is to count other subgraphs instead of edges. Given graphs $H$ and $G$, we let $\cN(H,G)$ denote the number of copies of $H$ in $G$. The \textit{generalized Tur\'an number} $\ex(n,H,F)$ is the largest $\cN(H,G)$ among $n$-vertex $F$-free graphs $G$.

Let $T_p(n)$ denote the {\it Tur\'an graph}; a balanced complete $p$-partite graph on $n$ vertices.  Tur\'an \cite{turan} proved that $T_p(n)$ is the unique extremal graph of $ex(n,K_{p+1})$, which is regarded as the beginning of the extremal graph theory.
The famous Erd\H{o}s-Stone-Simonovits Theorem \cite{Er,stone} states if $H$ is a graph with chromatic number $\chi(H)=\chi\ge 3$, then
$$\ex(n,H)=\left(\frac{\chi-2}{\chi-1}+o(1)\right){n\choose 2}.$$
That is,  the Tur\'an number $ex(n,H)$ is determined asymptotically for any nonbipartite graph $H$. However, it is still a challenging problem to  determine the exact value of the Tur\'an number and the extremal graphs for many nonbipartite graphs.

The \textit{friendship graph} or \textit{$k$-fan} $F_k$ consists of $k$ triangles all intersecting in one common vertex $v$.  Obviously, $F_k$ is nonbipartite. Erd\H{o}s, F\"uredi, Gould and Gunderson determined the Tur\'an number of it.
\begin{Theorem}(Erd\H{o}s, F\"uredi, Gould and Gunderson\cite{friendship})\label{friendship}
For every $k\ge 1$ and $n\ge 50k^2$,
$$ex\left(n,F_k\right)=\left\lfloor\frac{n^2}{4}\right\rfloor+\begin{cases}
k^2-k ~~~~{\it if } ~k~{\it is~ odd,}\\
k^2-\frac{3}{2}k~~{\it if } ~k~{\it is~ even.}
\end{cases}$$
\end{Theorem}

Recently, the problem of estimating generalized Tur\'{a}n number has received a lot of attention, some classical results have been extended to the generalized Tur\'an problem.
One can find them in \cite{Bollobas, chase,Gerbner,Grzesik,Hatami,Luo,Ma,Zykov}. A particular line of research is to determine for a given graph $H$, what graphs $F$ have the property that $\ex(n,H,F)=O(n)$. This was started by Alon and Shikhelman \cite{Alon}, who dealt with the case $H=K_3$, and was continued for other graphs in \cite{gerbner2,GP2017}.

An \textit{extended friendship graph} consists of $F_k$ for some $k\ge 0$ and any number of additional vertices or edges that do not create any additional cycles. Alon and Shikhelman \cite{Alon} showed that $\ex(n,K_3,F)=O(n)$ if and only if $F$ is an extended friendship graph. We remark that known results easily imply that if $F$ is not an extended friendship graph, then $\ex(n,K_3,F)=\omega(n)$ and it is also easy to see that adding further edges to $F$ without creating any cycle does not change linearity of $\ex(n,K_3,F)$. Hence the key part of their proof is the following theorem.

\begin{Theorem}(Alon and Shikhelman \cite{Alon}) For any $k$ we have $\ex(n,K_3,F_k)< (9k-15)(k+1)n$.
\end{Theorem}

 This upper bound for $ex(n,K_3,F_k)$ is not tight. For instance, for $k=2$, it was observed by Liu and Wang \cite{LW} that a hypergraph Tur\'an theorem of Erd\H os and S\'os \cite{Sos} gives the exact result for $ex(n,K_3,F_2)$. Let $\cF_k$ denote the 3-uniform hypergraph ($k$-star) consisting of $k$ hyperedges sharing exactly one vertex. Let $\ex_3(n,\cF_k)$ denote the largest number of hyperedges that an $\cF_k$-free $n$-vertex 3-uniform hypergraph can contain.
 
\begin{Theorem}(Erd\H{o}s and S\'os \cite{Sos}) For all $n\geq 3$,
\begin{displaymath}
\ex_3(n,\cF_2)=
\left\{ \begin{array}{l l}
n & \textrm{if\/ $n=4m$},\\
n-1 & \textrm{if\/ $n=4m+1$},\\
n-2 & \textrm{if\/ $n=4m+2$ or $n=4m+3$}.\\
\end{array}
\right.
\end{displaymath}
\end{Theorem}

Hence, it is interesting to determine the exact value of $ex(n,K_3,F_k)$ for any $F_k$  ($k\geq 3$).

Let $G=(V(G),E(G))$ be a connected simple graph and $e(G)=|E(G)|$. For any vertex $v\in V(G)$ and subset $S\subseteq V(G)$, let $N_S(v)$  denote the neighbors of $v$ in $S$ and $d_S(v)=|N_S(v)|$. If $S=V(G)$, then $N(v)=N_S(v)$ and $d(v)=d_S(v)$.  For $X,Y\subseteq V(G)$, $[X,Y]$ denotes the set of edges with one end in $X$ and another in $Y$ and $[x,Y]=[X,Y]$ if $X=\{x\}$. Let $\pi(G)$  denote the degree sequence of $G$. For two graphs $G_1$ and $G_2$, $G_1\cup G_2$ is the vertex disjoint union of $G_1$ and $G_2$ and $kG$ consists of $k$ copies of vertex disjoint union of $G$,  $G_1+G_2$ is the graph obtained by  taking $G_1\cup G_2$ and joining all pairs $v_1,v_2$ with $v_1\in V(G_1)$ and $v_2\in V(G_2)$. Let  $K_n$ and $\bar K_n$ denote  the complete graph and the empty graph on $n$ vertices, respectively.

We first define two graphs. Let $k\geq 4$ be even, $X=\{x_1,...,x_{k-1}\}$ and $Y=\{y_1,...,y_{k-1}\}$. The
graph $H_k'$ is a graph obtained from a complete bipartite graph with vertex classes $X$ and $Y$. We subdivide the edge $x_iy_i$ once for $i\leq \frac{k}{2}-1$, and then identify the $\frac{k}{2}-1$ inserted vertices into one vertex $z$. The graph $H_k$ is the complement of $H_k'$ deleting the edge $zy_{k/2}$, which is  shown in Figure 1.

\begin{picture}(70,90)(-76.5,75)
\centering
{\includegraphics[width=3.1in]{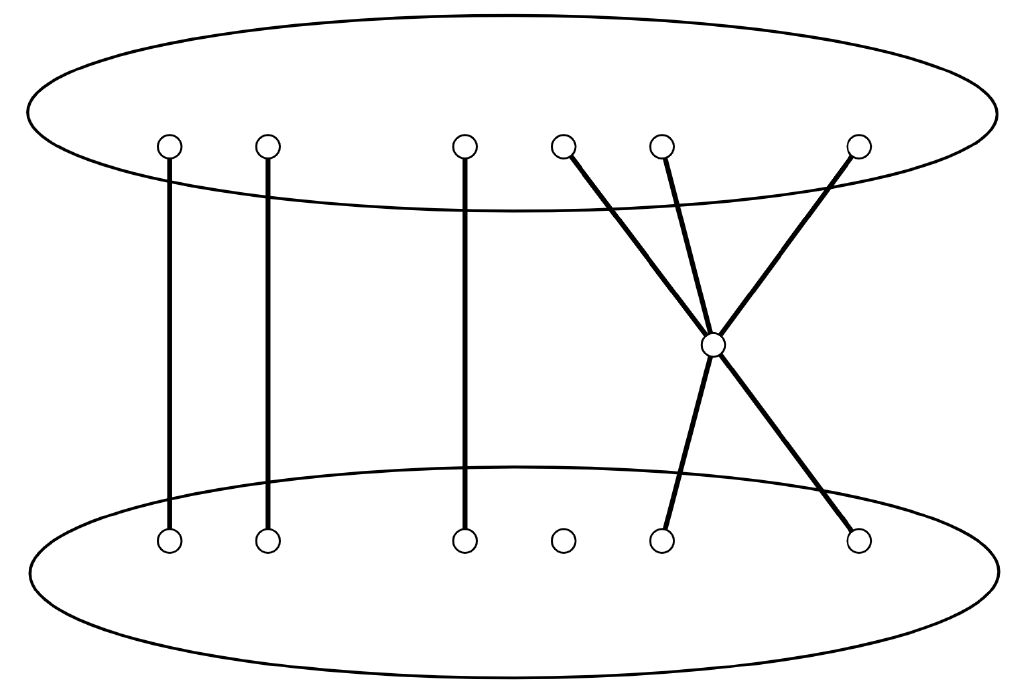}}
\put(-187,129){\makebox(3,3){$x_1$}}
\put(-135,134){\makebox(3,3){{\large$K_{k-1}$}}}

\put(-143,120){\makebox(3,3){$\cdots$}}
\put(-58,120){\makebox(3,3){$\cdots$}}

\put(-143,32){\makebox(3,3){$\cdots$}}
\put(-58,32){\makebox(3,3){$\cdots$}}

\put(-101,129){\makebox(3,3){$x_{k/2}$}}
\put(-36,129){\makebox(3,3){$x_{k-1}$}}

\put(-187,21){\makebox(3,3){$y_1$}}
\put(-135,15){\makebox(3,3){{\large$K_{k-1}$}}}

\put(-101,21){\makebox(3,3){$y_{k/2}$}}
\put(-36,21){\makebox(3,3){$y_{k-1}$}}

\put(-59,76){\makebox(3,3){$z$}}

\put(-113,-20){\makebox(3,3){Figure 1. The graph $H_k$}}
\end{picture}


\vskip 40mm

It is clear that $|H_k|=|H_k'|=2k-1$ and $\pi(H_k)=\pi(H_k')=(k-1,...,k-1,k-2)$.
\vskip 3mm
The main result of this paper is the following.

\begin{Theorem}\label{main}
Let $k\geq 3$ be an integer and $n\ge 4k^3$. If $k$ is odd, then
$$ex\big(n,K_3,F_k\big)=(n-2k)k(k-1)+2\binom{k}{3},$$
and $\bar K_{n-2k}+2K_k$ is the unique extremal graph, and if $k$ is even, then

$$ex\big(n,K_3,F_k\big)=(n-2k+1)k\left(k-\frac{3}{2}\right)+2\binom{k-1}{3}+\left(\frac{k}{2}-1\right)^2,$$
and $\bar K_{n-2k+1}+H_k$ is the unique extremal graph.
\end{Theorem}

\vskip 2mm

Given a graph $G$, we let $\cT(G)$ denote the 3-uniform hypergraph on the vertex set $V(G)$ where $\{u,v,w\}$ form a hyperedge if and only if $uvw$ is a triangle in $G$. The key observation is that if $G$ is $F_k$-free, then $\cT(G)$ is $\cF_k$-free. Therefore, $\ex(n,K_3,F_k)\le \ex_3(n,\cF_k)$. In the case $k=2$, the upper bound obtained this way matches the lower bound provided by $\lfloor n/4\rfloor$ vertex-disjoint copies of $K_4$, and in the case $n=4m+3$ we also have a triangle on the remaining vertices. This gives the exact value of $\ex(n,K_3,F_2)$.

The result of Erd\H os and S\'os \cite{Sos} was extended to arbitrary $k$ by Chung and Frankl \cite{Frankl}, after partial results \cite{Chung,ce,de}.

\begin{Theorem}[Chung and Frankl \cite{Frankl}]\label{chung} 
Let $k\ge 3$. If $n$ is sufficiently large, then
\begin{displaymath}
\ex_3(n,\cF_k)=
\left\{ \begin{array}{l l}
(n-2k)k(k-1)+2\binom{k}{3} & \textrm{if\/ $k$ is odd},\\
(n-2k+1)\frac{(2k-1)(k-1)-1}{2}+(2k-2)\binom{k-1}{2}+\binom{k-2}{2}-\frac{(k-2)(k-4)}{2}+\frac{k}{2} & \textrm{if\/ $k$ is even}.\\
\end{array}
\right.
\end{displaymath} 
and $\mathscr{F}_k:=\cT\big(\bar K_{n-2k}+2K_k\big)$ is the unique extremal 3-uniform hypergraph, when $k$ is odd.
\end{Theorem}

For odd $k$, this completes the proof of the upper bound. However, for even $k$, the construction giving the lower bound in the above theorem is not $\cT(G)$ for some $F_k$-free graph $G$. Still, the upper bound differs from the lower bound only by an additive constant $c(k)$. 
We will heavily use the tools provided by Chung and Frankl \cite{Frankl} to obtain the improvement needed in Theorem \ref{main}.


The rest of this paper is organized as follows. In Section 2, we present some preliminary results. In Section 3, we study the local structure within the neighborhood of a vertex in an  $F_k$-free graph, and some properties of a weight function defined on the vertices of triangles which is our main method for counting the number of triangles. These results  can be used to prove Theorem \ref{main} with the best coefficient of $n$ but a weak constant $f(k)$. The very technical Section 4 is devoted to  prove Theorem \ref{main} precisely. In Section 5, we give some concluding remarks. 

\section{Preliminaries}
As a preparation for proving our result, we first present some known theorems, and then we  count the number of triangles in a graph with given degree sequence, which are interesting of their own right.
\vskip 2mm
Let $\nu(G)$ denote the number of edges of a maximum matching in a graph $G$. The following is the famous result about the maximum matching,
\begin{Theorem}(Berge \cite{Berge})\label{Berge} Let $o(G-X)$ denote the number of odd components of $G-X$, then
$$\nu (G)=\frac{1}{2}\mathop{min}\Big\{|G|-o(G-X)+|X|:~X\subseteq V(G)\Big\}.$$
\end{Theorem}



\begin{Theorem}(Chung and Frankl \cite{Frankl})\label{Frankl}
Let $k$ be an even integer and $H$ be a graph on $2k-1$ vertices and with $\pi(H)=(k-1,\ldots,k-1,k-2)$, then either
$$\cN(K_3,H)\ge \left(\frac{k}{2}-1\right)^2-1,$$
or
$$\cN(K_3,H)=\left(\frac{k}{2}-2\right)\left(\frac{k}{2}-1\right) ~\text{and}~ H=H'_k.$$
\end{Theorem}

\begin{Theorem}\label{thm3}
Let $k$ be an even integer and $H$ be a graph on $2k-1$ vertices with $\pi(H)=(k-1,\ldots,k-1,k-2)$, then
$$\cN(K_3,H)\le 2\binom{k-1}{3}+\left(\frac{k}{2}-1\right)^2,$$
equality holds if and only if $H=H_k$.
\end{Theorem}

\pf The proof will be similar to the proof of Goodman. It is easy to see that
 $$\mathcal{N}(K_3,H_k)=2\binom{k-1}{3}+\binom{k/2}{2}+\binom{k/2-1}{2}=2\binom{k-1}{3}+\left(\frac{k}{2}-1\right)^2.$$
Let $\overline{H}$ be the complement of $H$. For any triple $(x,y,z)$, if  $xyz$ is neither a triangle in $H$ nor a triangle in $\overline{H}$, then it is easy to check exactly two of the three, say $x,y$ such that $|[x,\{y,z\}]|=1$ and $|[y,\{x,z\}]|=1$ in $H$. Thus, we have
\[\begin{split}
\mathcal{N}(K_3,H)&=\binom{2k-1}{3}-\mathcal{N}\!\left(K_3,\overline{H}\right)-\frac{1}{2}\sum_v d(v)(2k-2-d(v))\\
&=\binom{2k-1}{3}-(k-1)^3-\frac{1}{2}k(k-2)-\mathcal{N}\!\left(K_3,\overline{H}\right)\\
&=2\binom{k-1}{3}+\frac{(k-2)^2}{2}\mathcal{N}\!\left(K_3,\overline{H}\right).
\end{split}\]
Obviously, it is sufficient  to show $\mathcal{N}\!\left(K_3,\overline{H}\right)\ge (\frac{k}{2}-1)^2$.

Note that $\overline{H}$ is a graph on $2k-1$ vertices with $\pi\!\left(\overline{H}\right)=(k-1,\ldots,k-1,k)$. Let $z$ be the vertex of degree $k$ in $\overline{H}$.

If there is an edge $zz'\in E\!\left(\overline{H}\right)$ such that $zz'$ is contained in at least two triangles, then $\mathcal{N}\!\left(K_3,\overline{H}\right)\ge \mathcal{N}\!\left(K_3,\overline{H}-zz'\right)+2$. Because of $\pi\!\left(\overline{H}-zz'\right)=(k-1,\ldots,k-1,k-2)$, by Theorem \ref{Frankl}, either $\mathcal{N}(K_3,\overline{H}-zz')\ge \left(\frac{k}{2}-1\right)^2-1$ or $\overline{H}-zz'=H'_k$.
In the former case, we have
$$\mathcal{N}(K_3,H)\le 2\binom{k-1}{3}+\left(\frac{k}{2}-1\right)^2-1.$$
In the latter case, it is easy to check that $\overline{H}=\overline{H_k}$, and hence $H=H_k$.

If each edge $zz'$ in $E\!\left(\overline{H}\right)$ is contained in at most one triangle, then $\Delta\!\left(\overline{H}[N(z)]\right)\le 1$. Let $V_1=N(z)$, $V_2=V(H)-V_1$, $s=e\!\left(\overline{H}[V_1]\right)$ and $t=e\!\left(\overline{H}[V_2]\right)$. Count the edges between $V_1$ and $V_2$ in two ways, we have
$$k(k-1)-2s=(k-1)^2+1-2t,$$  which implies $s-t=\frac{k}{2}-1$. However, because $s\le \frac{k}{2}$, we must have $s=\frac{k}{2}$ and $t=1$, or $s=\frac{k}{2}-1$ and $t=0$. Since each edge in $\overline{H}[V_1]$ can form a triangle with at least $k-3$ vertices in $V_2$, we get
$$\mathcal{N}\!\left(K_3,\overline{H}\right)\ge  \frac{k}{2}(k-3)+(k-4)>\left(\frac{k}{2}-1\right)^2$$
in the former case,  and
$$\mathcal{N}\!\left(K_3,\overline{H}\right)\ge \left(\frac{k}{2}-1\right)(k-3)\ge \left(\frac{k}{2}-1\right)^2$$
in the latter case with equality only if $k=4$.  In this case,  it is not difficult to check that $\overline{H}=\overline{H_4}$, and so $H=H_4$. $\hfill\blacksquare$

\begin{Theorem}\label{thm4}
Let $k$ be an even integer and $H$ be a graph on $2k-1-2s$ vertices with  $\pi(H)=(k-1,\ldots,k-1,k-2)$, then
$$\cN(K_3,H)\le \frac{1}{6}(2k-1-2s)\Big((k-1)(k-2)-(k-1-2s)(2s+1)\Big)+\frac{1}{2}-s.$$
\end{Theorem}
\pf Let $\varLambda(H)$ denote the number of triples $(x,y,z)$ having exactly two edges in $H$, say  $xy,xz\in E(H)$  and $yz\notin E(H)$. Because $|N(y)\cap N(z)|\ge d(y)+d(z)-(|H|-2)$ for every nonadjacent pair $(y,z)$,  and there are $|H|-(d(y)+1)$ nonadjacent pairs containing $y$ for any $y\in V(H)$, we have
$$\varLambda(H)\ge \frac{1}{2}\Big((2k-1-2s)(k-1-2s)+1\Big)(2s+1)-(k-2s).$$
 On the other hand, since
 $$(2k-2-2s)\binom{k-1}{2}+\binom{k-2}{2}=\varLambda(H)+3\cN(K_3,H),$$
we get
$$3\cN(K_3,H) \le \frac{1}{2}(2k-1-2s)\Big((k-1)(k-2)-(k-1-2s)(2s+1)\Big)+\frac{3}{2}-3s. $$
This completes the proof. $\hfill\blacksquare$

\section{Some properties of $F_k$-free graphs and a weight function}

Let  $G$ be an $F_k$-free graph, $uv\in E(G)$ and $N(uv)=N(u)\cap N(v)$. Clearly, $|N(uv)|$ is the number of triangles containing the edge $uv$ in $G$. We classify the edges into the following three classes:
\vskip 2mm
$\bullet$ Heavy edges: $\mathcal{H}~=\big\{uv: |N(uv)|\ge 2k-1\big\}$,

$\bullet$ Medium edges: $\mathcal{M}=\big\{uv:k\le|N(uv)|\le 2k-2 \big\}$, and 

$\bullet$ Light edges: $\mathcal{L}~=\big\{uv:1\le |N(uv)|\le k-1\big\}$.
\vskip 2mm
\noindent For a fixed vertex $u\in V(G)$, let $G_u=G[N(u)]$ and
\vskip 2mm
$\bullet$ ~\!$\mathcal{H}(u)~=\big\{v:v\in N(u)~and~uv\in \mathcal{H}\big\}$.

$\bullet$ $\mathcal{M}(u)=\big\{v:v\in N(u)~and~uv\in \mathcal{M}\big\}$,

$\bullet$ ~\!$\mathcal{L}(u)~=\big\{v:v\in N(u)~and~uv\in \mathcal{L}\big\}$,
\vskip 2mm
\noindent This notation will be used throughout the rest of this paper.

Since $G$ is $F_k$-free, then $\nu(G_u)\le k-1$ for any $u$. Thus, Theorem \ref{Berge} implies
\begin{Observation}
There exists some $X\subseteq V(G_u)$ such that
\begin{equation}\label{eq3.1}
\sum\limits_{i=1}^\ell\left\lfloor \frac{|C_i|}{2}\right\rfloor+|X|\le k-1,
\end{equation}
where $C_1,...,C_{\ell}$ are all the components of $G_u-X$.
\end{Observation}

\begin{Lemma}\label{lemma1}\label{L1}  Let $G$ be an $F_k$-free graph, $u\in V(G)$ and $X$ a subset of $V(G_u)$ satisfying the equation (\ref{eq3.1}). Then we have the following:
\vskip 2mm
\noindent (i)~~$\mathcal{H}(u)\subseteq X$. Moreover, $|\mathcal{H}(u)|\le k-1$ and if equality holds, then $\mathcal{M}(u)=\emptyset$.

\noindent (ii)~$|\mathcal{H}(u)|+\frac{1}{2}|\mathcal{M}(u)|\le k-\frac{1}{2}$.
\end{Lemma}
\pf Let $C_1,...,C_{\ell}$ be the components of $G_u-X$.
\vskip 2mm
(i) Let $v\in \mathcal{H}(u)$, we know that $|N(uv)|\ge 2k-1$. If $v$ lies in some component $C_i$, then $N(uv)\subseteq V(C_i)\cup X$ and so
$$\frac{1}{2}\big|(C_i\cap N(uv))\cup \{v\}\big|+\big|X\cap N(uv)\big|\ge k,$$
which contradicts (\ref{eq3.1}). Hence we have $\mathcal{H}(u)\subseteq X$.

By (\ref{eq3.1}), we have $|\mathcal{H}(u)|\le k-1$, and if $|\mathcal{H}(u)|=k-1$, then $X=\mathcal{H}(u)$ and  $|C_i|=1$ for $1\leq i\leq \ell$.  Let $v$ be any vertex of $G_u-X$, then $N(uv)\subseteq X$ and hence $uv\in \mathcal{L}$, and so $\mathcal{M}(u)=\emptyset$.
\vskip 2mm
(ii) Clearly,  $N(uv)\subseteq X\cup V(C_i)$ if $v\in V(C_i)$. Thus, if there are two components, say $C_1,C_2$, such that $\mathcal{M}(u)\cap V(C_i)\not=\emptyset$, then
$|X|+|C_i|\ge k+1$ for $i=1,2$. Hence we have   $|X|+\lfloor |C_1|/2\rfloor+\lfloor |C_2|/2\rfloor\ge k$, which contradicts (\ref{eq3.1}).
Thus, we may assume $\mathcal{M}(u)\subseteq X\cup V(C_i)$. Note that $\mathcal{H}(u)\subseteq X$ as shown in (i), and
\[\begin{split}
\big|\mathcal{H}(u)\big|+\frac{1}{2}\big|\mathcal{M}(u)\big|&=\big|\mathcal{H}(u)\big|+\frac{1}{2}\big|\mathcal{M}(u)\cap (X-\mathcal{H}(u))\big|+\frac{1}{2}\big|\mathcal{M}(u)\cap C_i\big|\\
&\le \big|X\big|+\frac{1}{2}\big|C_i\big|\le \big|X\big|+\left\lfloor \frac{|C_i|}{2}\right\rfloor+\frac{1}{2}\leq k-\frac{1}{2}.
\end{split}\]
The proof of the lemma is complete.$\hfill\blacksquare$

\vskip 5mm
For each triangle $T=xyz$ in $G$,
assign $T$ of weight $1$ and define a distribution rule $w(T,\cdot)$ to distribute the weight 1 to its three vertices as below(suppose $|N(xy)|\ge |N(yz)|\ge |N(xz)|$):
\vskip 2mm
\hskip 5mm $w(T,x)=w(T,y)=w(T,z)=\frac{1}{3},~~~~~~{\it if} ~E(T)\cap \mathcal{H}=\emptyset~ {\it or} ~E(T)\cap \mathcal{L}=\emptyset,$
\vskip 1mm
\hskip 5mm $w(T,x)=w(T,z)=\frac{1}{2},~ w(T,y)=0,\,~~{\it if}~xy\in \mathcal{H},~yz\in \mathcal{H}\cup \mathcal{M}~\text{and}~xz\in \mathcal{L},$
\vskip 1mm
\hskip 5mm $w(T,x)=w(T,y)=0,~w(T,z)=1,~~~{\it if}~xy\in \mathcal{H} ~\text{and}~yz,xz\in \mathcal{L}.$
\vskip 2mm

\noindent
Now, we define a weight function $f(u)$ for each vertex $u$ of $G$ as follows.
$$f(u)=\sum\limits_{vx\in E(G_u)}w(uvx,u)$$
if $u$ lies in triangles, and $f(u)=0$ otherwise.
It is clear
$$\cN(K_3,G)=\sum_{u\in V(G)}f(u).$$
Now,  we first discuss some properties of the weight functions $w(T,\cdot)$ and $f(u)$.

\begin{Lemma} \label{lemma2}Let $uv$ be an edge of an $F_k$-free graph with $k\geq 3$, then either
$$\sum_{x\in N(uv)}w(uvx,u)=k-1,$$
 if $uv\in \mathcal{L}, |N(uv)|=k-1, vx\in \mathcal{H}$ and $ux\in \mathcal{L}$ for any $x\in N(uv)$,  or

$$\sum_{x\in N(uv)}w(uvx,u)\leq k-\frac{3}{2}
$$
otherwise.


\end{Lemma}
\pf Let $\mathcal{H}'(v)$, $\mathcal{M}'(v)$ and $\mathcal{L}'(v)$ be  $\mathcal{H}(v)$, $\mathcal{M}(v)$ and $\mathcal{L}(v)$ intersecting with $N(uv)$, respectively.
It is clear
$$\sum_{x\in N(uv)}w(uvx,u)=\sum_{x\in \mathcal{H}'(v)}w(uvx,u)+\sum_{x\in \mathcal{M}'(v)}w(uvx,u)+\sum_{x\in \mathcal{L}'(v)}w(uvx,u).$$

We distinguish three cases on the number of $|N(uv)|$.
\vskip 2mm
\noindent {\bf Case 1.} $uv\in \mathcal{H}$
\vskip 2mm
By the definition of $w(T,\cdot)$, $w(uvx,u)\leq \frac{1}{2}$ if $x\in \mathcal{H}'(v)\cup \mathcal{M}'(v)$ and $w(uvx,u)=0$ if $x\in \mathcal{L}'(v)$.
Noting that  $u\in \mathcal{H}(v)-\mathcal{H}'(v)$, we have $|\mathcal{H}'(v)|+ \frac{1}{2}|\mathcal{M}'(v)|\le k-\frac{3}{2}$ by  Lemma \ref{L1}(ii), and hence
\begin{equation}\label{eq3.2}
\sum_{x\in N(uv)}w(uvx,u)\leq \frac{1}{2}\big|\mathcal{H}'(v)\big|+\frac{1}{2}\big|\mathcal{M}'(v)\big|\le k-\frac{3}{2}-\frac{1}{2}\big|\mathcal{H}'(v)\big|,
\end{equation}
which implies the result holds.
\vskip 2mm
\noindent {\bf Case 2.}  $uv\in \mathcal{M}$
\vskip 2mm
In this case, $|N(uv)|\le 2k-2$. Since $uv\in \mathcal{M}$ implies $\mathcal{M}(v)\not=\emptyset$, by Lemma \ref{L1}(i), we have $|\mathcal{H}'(v)|\le |\mathcal{H}(v)|\le k-2 $.  By the definition of $w(T,\cdot)$, $w(uvx,u)\leq \frac{1}{2}$ if $x\in \mathcal{H}'(v)$ and $w(uvx,u)\leq \frac{1}{3}$ if $x\in \mathcal{M}'(v)\cup \mathcal{L}'(v)$.
Thus, we have

\begin{equation}\label{eq3.3}
\sum_{x\in N(uv)}w(uvx,u)\le \frac{1}{2}\big|\mathcal{H}'(v)\big|+\frac{1}{3}\big|\mathcal{M}'(v)\big|+\frac{1}{3}\big|\mathcal{L}'(v)\big|\\
\le k-1-\frac{k}{6}.
\end{equation}
The upper bound $k-\frac{3}{2}$ follows from the assumption $k\geq 3$.
\vskip 2mm
\noindent {\bf Case 3.} $uv\in \mathcal{L}$
\vskip 2mm
Because $uv\in \mathcal{L}$, we have $|N(uv)|\le k-1$. By the definition of $w(T,\cdot)$, some triangles satisfy $w(uvx,u)=1$ and other triangles satisfy $w(uvx,u)\leq \frac{1}{2}$. Thus we have
\begin{equation}\label{eq3.4}
\sum_{x\in N(uv)}w(uvx,u)\le (k-1)-\frac{1}{2} \left|\left\{uvx: w(uvx,u)\le \frac{1}{2}\right\}\right|,
\end{equation}
which implies $\sum_{x\in N(uv)}w(uvx,u)\leq k-\frac{3}{2}$  or $\sum_{x\in N(uv)}w(uvx,u)=k-1$, and the latter holds if and only if $|N(uv)|=k-1$, and all triangles $uvx$ satisfy $w(uvx,u)=1$, that is, $vx\in \mathcal{H}$ and $ux\in \mathcal{L}$ for any $x\in N(uv)$.
 $\hfill\blacksquare$

\begin{Lemma}\label{lem3}  Suppose $G$ is an $F_k$-free graph and $k\geq 4$ is even. Let $u\in V(G)$, $X$ be a subset of $V(G_u)$ satisfying (\ref{eq3.1}) and $C_1,...,C_\ell$ be the components of $G_u-X$ with $|C_1|\geq \cdots\geq |C_\ell|$. Then

$$f(u)\le k\left(k-\frac{3}{2}\right)-\frac{1}{2},$$
or $f(u)=k\left(k-\frac{3}{2}\right)$ and the following hold:
\vskip 2mm

(i) If $X\not=\emptyset$, then $X$ is an independent set and $d_{G_u}(v)=k-1$ for any $v\in X$;

(ii) $\pi(C_1)=(k-1,\ldots,k-1,k-2)$, and either $G_u=C_1\cup K_{k-1}$ with $|C_1|=k+1$, or
$G_u-X=C_1\cup (\ell-1)K_1$ with $|C_1|=2k-1-2|X|\geq k+1$;

(iii) $E(G_u)\subseteq \mathcal{H}$, $[u, G_u]\subseteq \mathcal{L}$ and $\Delta(G_u)=k-1$.
\end{Lemma}
\noindent{\bf Remark.} There is a similar lemma in Chung and Frankl's paper\cite{Frankl} when they deal with function $\ex_3\big(n,\mathcal{F}_k\big)$. However, in their lemma, they overlooked the case $G_u=C_1\cup K_{k-1}$.
Using our method in Section 4, it is not difficult to complete the proof of this missed case, too. 

\pf By  (\ref{eq3.1}), $|C_i|\leq k$ for all $i\not=1$. Let $uvx$ be any triangle. Then

$$f(u)=\sum_{i=1}^\ell\sum_{vx\subseteq E(C_i)}w(uvx,u)+\sum_{\{v,x\}\cap X\neq \emptyset}w(uvx,u).$$

\vskip 3mm
If the edge $vx$ satisfies $\{v,x\}\cap X\not=\emptyset$, then
by Lemma \ref{lemma2}, we have

\begin{equation}\label{eq3.5}
\sum_{\{v,x\}\cap X\neq \emptyset}w(uvx,u)\le \sum_{v\in X}\sum_{x\in N(uv)}w(uvx,u)\leq \big|X\big|(k-1).
\end{equation}

\vskip 3mm
If $vx\in E(C_i)$ with $|C_i|\leq k$, then noting that $k$ is even, $uvx$ is a triangle for each $vx\in E(C_i)$ and $w(uvx,u)\leq 1$, we have
$$\sum_{vx\in E(C_i)}w(uvx,u)\le \frac{1}{2}\sum_{v\in V(C_i)}\sum_{x\in N(uv)\cap C_i} w(uvx,u)~~~~~~~~~~~~~~~~~~~$$
\begin{equation}\label{eq3.6}~~~~~~~~~~\le \frac{1}{2}\big|C_i\big|\left(\big|C_i\big|-1\right)\leq  \left\lfloor\frac{|C_i|}{2}\right\rfloor(k-1).
\end{equation}

If $|C_1|\leq k$, that is, $|C_i|\leq k$ for $1\leq i\leq \ell$, then by (\ref{eq3.1}), (\ref{eq3.5}) and (\ref{eq3.6}), we have
$$f(u)\le (k-1)\left(\sum_{i=1}^\ell\left\lfloor \frac{|C_i|}{2}\right \rfloor+\big|X\big|\right)\le (k-1)^2<k\left(k-\frac{3}{2}\right)-\frac{1}{2}.$$
The last inequality holds because of $k\geq 4$. So, we may assume that $|C_1|>k$. 

If $\Delta(C_1)\geq k$, say $d_{C_1}(v)\ge k$ for some vertex $v$ in $C_1$, then $|N(uv)|\geq k$ and so $uv\notin \mathcal{L}$. Since $\mathcal{H}(u)\subseteq X$ by Lemma \ref{L1}(i), we have $v\notin \mathcal{H}(u)$ and hence $uv\in \mathcal{M}$.  By (\ref{eq3.3}), we have
$\sum_{x\in N(uv)\cap C_1} w(uvx,u)\le (k-1)-\frac{k}{6}.$
Meanwhile, because $d_{C_1}(v)\geq k\geq 4$, there exists $v_i\in N(uv)\cap C_1$ for $1\leq i\leq 4$. Note that $v\in \mathcal{M}(u)$ and $v\in N(uv_i)$,  by  Lemma \ref{lemma2}, we have
$\sum_{x\in N(uv_i)\cap C_1} w(uv_ix,u)\le (k-1)-\frac{1}{2}~\text{for}~1\leq i\leq 4$, and
\begin{equation}\label{eq3.7}
\sum_{vx\in E(C_1)}w(uvx,u)=\frac{1}{2}\sum_{v\in C_1}\sum_{x\in N(uv)\cap C_1}w(uvx,u)\leq \frac{1}{2}\big|C_1\big|(k-1)-\left(\frac{k}{12}+1\right).
\end{equation}
If $\Delta(C_1)\le k-1$, then
\begin{equation}\label{eq3.8}
\sum_{vx\subseteq E(C_1)}w(uvx,u)=\frac{1}{2}\sum_{v\in C_1}\sum_{x\in N(uv)\cap C_1}w(uvx,u) \le \left\lfloor\frac{1}{2}\big|C_1\big|(k-1)\right\rfloor.
\end{equation}
Set $\mu(C_1)=\frac{k}{12}+1$ if $\Delta(C_1)\geq k$, $\mu(C_1)=\frac{1}{2}$ if $|C_1|$ is odd and $\Delta(C_1)\leq k-1$ and $\mu(C_1)=0$ if $|C_1|$ is even and $\Delta(C_1)\leq k-1$, then   (\ref{eq3.5})-(\ref{eq3.8}) imply
\begin{equation}\label{eq3.9}
f(u)\le (k-1)\left(\frac{|C_1|}{2}+\sum_{i=2}^\ell\left\lfloor \frac{|C_i|}{2} \right\rfloor+\big|X\big|\right)-\mu\big(C_1\big).
\end{equation}

Assume that $f(u)>k\left(k-\frac{3}{2}\right)-\frac{1}{2}$. By (\ref{eq3.1}) and (\ref{eq3.9}), we have $\mu(C_1)=\frac{1}{2}$.  In this case, $|C_1|\geq k+1$ is odd and $\Delta(C_1)\leq k-1$.
Note that if one of the equalities  in (\ref{eq3.5}), (\ref{eq3.6}) and (\ref{eq3.8}) does not hold, then the upper bound in (\ref{eq3.9}) can be reduced by at least an extra $\frac{1}{2}$. This implies the equalities  in (\ref{eq3.5}), (\ref{eq3.6}) and (\ref{eq3.8}) holds.

It is clear that the equalities  in (\ref{eq3.5}) hold if and only if $X$ is an independent set and  $\sum_{x\in N(uv)}w(uvx,u)=k-1$ for any $v\in X$. By Lemma \ref{lemma2}, we get that $d_{G_u}(v)=k-1$, $uv,ux\in \mathcal{L}$ and $vx\in \mathcal{H}$ for any $v\in X$.

Since $|C_1|\geq k+1$ is odd, $\Delta(C_1)\leq k-1$ and the equality in (\ref{eq3.8}) holds, we can deduce that  $\pi(C_1)=(k-1,k-1,\ldots,k-2)$, $E(C_1)\subseteq\mathcal{H}$ and $[u,C_1]\subseteq \mathcal{L}$.

Because  equality (\ref{eq3.6}) holds, recalling $|C_1|\geq k+1$ and $k$ is even, by (\ref{eq3.1}), we have $|C_i|\in\{1,k-1\}$ for $i\ge 2$ and each $C_i$ is a clique with $E(C_i)\subseteq \mathcal{H}$ and $[u,C_i]\subseteq \mathcal{L}$. In addition, if  $|C_i|=k-1$ for some $i\ge 2$, then by (\ref{eq3.1}), $|C_1|=k+1$ and $X=\emptyset$, that is, $G_u-X=C_1\cup K_{k-1}\cup (\ell-2)K_1$. If $|C_i|=1$ for all $i\geq 2$, then $|C_1|=2k-1-2|X|$.

So, the statements  (i), (ii) and (iii) hold.$\hfill\blacksquare$
%
%
%

\vskip 5mm

\begin{definition}
For any vertex $u\in V(G)$, the loss of $u$ is the number
$$k\left(k-\frac{3}{2}\right)-f(u).$$
\end{definition}

See the following simple observations about the losses.

\begin{Observation}\label{O2}
If some vertex $v\in X$ has $\sum_{x\in N(uv)}w(uvx,u)\le(k-1)-c$, then the edge $uv$ contributes $c$ to the loss of $u$.
\end{Observation}

\pf It is a direct consequence of (\ref{eq3.5}). $\hfill\blacksquare$

\begin{Observation}\label{O3}
An edge $uv\in \mathcal{H}$ contributes $\frac{1}{2}$ to the loss of $u$. Moreover, a triangle $uvx$ with $uv,vx\in \mathcal{H}$  contributes another $\frac{1}{2}$ to the loss of $u$.
\end{Observation}
\pf Since $uv\in \mathcal{H}$, by  (\ref{eq3.2}), $\sum_{x\in N(uv)}w(uvx,u)\le k-\frac{3}{2}-\frac{1}{2}|\mathcal{H}'(v)|$. Because $v\in X$ by Lemma \ref{lemma1}, the edge $uv$ contributes $\frac{1}{2}$ to the loss of $u$ by Observation \ref{O2}. Moreover, since a triangle $uvx$ with $uv,vx\in \mathcal{H}$ satisfies $x\in \mathcal{H}'(v)$,  so it contributes another $\frac{1}{2}$ to the loss of $u$ by (\ref{eq3.2}).  $\hfill\blacksquare$

\begin{Observation}\label{O4}
Let  $uv\in \mathcal{M}$. If $\sum_{x\in N(uv)} w(uvx,u)\le (k-1)-c$, then the edge $uv$ contributes at least
$min\left\{\frac{c}{2}, ~\frac{k}{4}-\frac{1}{2}\right\}$
to the loss of $u$. Moreover, the edge $uv$ contributes at least $\frac{k}{12}$ to the loss of $u$.
\end{Observation}
\pf  If $v\in X$, then by (\ref{eq3.3}) and Observation \ref{O2}, $uv$ contributes $c$ to the loss of $u$.
If $v\in V(C_1)$ and $|C_1|\geq k+1$, then  by  (\ref{eq3.7}) and (\ref{eq3.8}), $uv$ contributes at least $\frac{c}{2}$ to the loss of $u$.
If $v\in V(C_i)$ for  some $i$ with $|C_i|\le k$, then by (\ref{eq3.6}), there is a gap between  $\left\lfloor \big|C_i\big|/2\right\rfloor(k-1)$ and $\frac{1}{2}\big|C_i\big|\big(\big|C_i\big|-1\big)$, and for this gap,  any edge $uv'$ with $v'\in V(C_i)$  contributes
$$\frac{1}{\big|C_i\big|}\left(\left\lfloor\frac{|C_i|}{2}\right\rfloor(k-1)-\frac{1}{2}\big|C_i\big|\left(\big|C_i\big|-1\right)\right)$$
to the loss of $u$. On the other hand, because
$$\sum_{x\in N(uv)\cap C_i} w(uvx,u)\le \frac{1}{2}\big(\big|C_i\big|-1\big)=\big(\big|C_i\big|-1\big)-\frac{1}{2}\big(\big|C_i\big|-1\big),$$
this reduces the right hand of (\ref{eq3.6}) by an additional $\frac{1}{4}\big(\big|C_i\big|-1\big)$. Hence the total loss of $u$ contributed by the edge $uv$ is at least
\[\begin{split}
&\frac{1}{\big|C_i\big|}\left(\left\lfloor\frac{|C_i|}{2}\right\rfloor(k-1)-\frac{1}{2}\big|C_i\big|\left(\big|C_i\big|-1\right)\right)+\frac{1}{4}\big(\big|C_i\big|-1\big)
\ge \frac{k}{4}-\frac{1}{2}.
\end{split}\]
Together with (\ref{eq3.3}), $c\ge \frac{k}{12}$, it implies that the statements of the lemma are proved. $\hfill\blacksquare$


\section{Proof of Theorem \ref{main}}

Let $G$ be an extremal graph of $ex(n,K_3,F_k)$.

If $k$ is odd, then by Theorem  \ref{chung}, we have
$$\mathcal{N}(K_3,G)= e\!\left(\cT(G)\right)\leq \ex_3(n,\cF_k)= (n-2k)k(k-1)+2\binom{k}{3},$$
and the unique extremal hypergraph is $\mathscr{F}_k$ for which equality holds. Because
$$\mathcal{N}(K_3, \bar K_{n-2k}+2K_k)= (n-2k)k(k-1)+2\binom{k}{3},$$
and $\cT(\bar K_{n-2k}+2K_k)=\mathscr{F}_k$, we get
$$\mathcal{N}(K_3, G)= (n-2k)k(k-1)+2\binom{k}{3},$$
where equality holds if and only if $G=\bar K_{n-2k}+2K_k$.

%
\vskip 3mm

The remaining part is devoted to the case when $k\geq 4$ is even. Because an edge not lying in a triangle makes no contribution to $\cN(K_3, G)$, we may assume that each edge of $G$ is covered by some triangles.

If $f(v)=k\left(k-\frac{3}{2}\right)$, then we call $v$ a {\it good} vertex. Let $U_1=\{v:v~is~ good\,\}$. Since $n\geq 4k^3$ and $f(v)\leq k\left(k-\frac{3}{2}\right)-\frac{1}{2}$ for any $v\notin U_1$ by Lemma \ref{lem3}, we have
\[\begin{split}
nk\left(k-\frac{3}{2}\right)-\frac{1}{2}\left(n-|U_1|\right) &\ge\sum_{v\in V(G)} f(v)=\mathcal{N}(K_3,G)\\
  &\ge (n-2k+1)k\left(k-\frac{3}{2}\right)+2\binom{k-1}{3}+\left(\frac{k}{2}-1\right)^2,
 \end{split}\]
which implies
$$|U_1|\ge n-2(2k-1)k\left(k-\frac{3}{2}\right)> 2k.$$
Moreover, if there exist $v,v'\in U_1$ such that $vv'\in E(G)$, then $vv'\in \mathcal{L}$ by Lemma \ref{lem3}. Let $vv'x$ be a triangle.  Applying Lemma \ref{lem3} to $v$, we have  $v'x\in \mathcal{H}$ and using Lemma \ref{lem3} to $v'$, we have $v'x\in \mathcal{L}$, a contradiction. Therefore, $U_1$ is an independent set.
\vskip 1mm
Let $u\in U_1$ be given, $G_u=G[N(u)]$ as before and $U_2=V(G)-V(G_u)-U_1$. We will prove Theorem \ref{main} by showing $G$ is an extremal graph only if $U_2=\emptyset$. Since the proof is a little complicated and long, so we sketch it first in the following two paragraphs.

In the case when $N(u')=N(u)$ for any $u'\in U_1-\{u\}$,
our main idea for doing this is to partition the total weights of all vertices of $G_u$ into two parts: One part comes from the triangles contained in $G_u$, which is exactly $\mathcal{N}(K_3, G_u)$, and another part is contributed by the triangles containing one or two vertices in $U_2$.  And then we use discharge method to transfer the latter part to the vertices in $U_2$ such that $f(v)\leq k\left(k-\frac{3}{2}\right)$ is still valid for any $v\in U_2$ after transferring. Using this method, we show that if $U_2\not=\emptyset$, then the total weight of $G$ is less than the expected number.

 In the case when  there is some $u'\in U_1-\{u\}$  such that $N(u)\not=N(u')$,  we  transform $G$ into a  graph $G'$  such that $G'$ and $G$ have the same good vertices, and all good vertices of $G'$ have the same neighborhood as $N(u)$, through an operation as follows: Delete all the edges between $u'$ and $G_{u'}$ and add  new edges joining $u'$ to all vertices in $G_u$. Repeat this operation until all  vertices in $U_1$ have the same neighborhood $N(u)$. Let $G'$ be the resulting graph, $U_1'=\{v:v~is~good~in~G'\}$ and $U_2'=V(G')-V(G'_u)-U_1'$. We will see that $\mathcal{N}(K_3, G')=\mathcal{N}(K_3, G)$, $G'$ is also $F_k$-free and $U_1'=U_1$.

Firstly, since $u'$ is good, by Lemma \ref{lem3}, we have $f(u')=e(G_{u'})=k\!\left(k-\frac{3}{2}\right)$, which implies we destroy $k\!\left(k-\frac{3}{2}\right)$ triangles first and then add $k\!\left(k-\frac{3}{2}\right)$ new triangles, and so  $\mathcal{N}(K_3, G')=\mathcal{N}(K_3, G)$. Moreover, $G'_u=G_u$. Secondly,
since $G_u$ has no $kK_2$ and $\Delta(G_u)=k-1$ by Lemma \ref{lem3}, we can see that $G'$ is also $F_k$-free after an easy check. Finally, because $|U_1|>2k$, we have $E(G'_u)\subseteq \mathcal{H}$, which implies $v\notin U_1'$ for any $v\in V(G'_u)$ by Lemma \ref{lem3}. Furthermore,  since $\Delta(G'_u)=k-1$ and $U_1'$ is an independent set, $[U_1,G'_u]\subseteq \mathcal{L}$. Thus, we have $U_1\subseteq U_1'$ by the definition of $w(T,\cdot)$. Suppose that there is some $v\in U_2$ in $G$ such that $v\in U_1'$ in $G'$. Let $G'_v=G'[N(v)]$ and $X'\subseteq G'_v$ satisfy (\ref{eq3.1}). Since $v\in U_1'$,  $E(G'_v)\subseteq \mathcal{H}$ and $[v,G'_v]\subseteq \mathcal{L}$ by Lemma \ref{lem3}. 
By the operation above, $E(G'_v)\subseteq \mathcal{H}$ in $G$. Since $v\in U_2$, there is some $v'\in G_v'$ such that $vv'\notin \mathcal{L}$ in $G$, which means there is some $u'\in U_1$ such that $u'vv'$ is a triangle in $G$. Note that $V(G_v')\cup \{u'\}\subseteq N_G(v)$. If $v'\notin X'$, then by Lemma \ref{lem3}, it is easy to check that $G[\{u'\}\cup V(G_v')]$ contains $kK_2$, and so $G$ has an $F_k$. Thus we have $v'\in X'$. In this case, by Lemma \ref{lem3}, $|\mathcal{H}(v')|=k-1$ in $G$. Let $X''\subseteq G_{v'}$ satisfy (\ref{eq3.1}). By Lemma \ref{lemma1}, $\mathcal{H}(v')\subseteq X''$ and hence $|X''|=k-1$. Because $u'v$ is an edge in $G_{v'}-X''$, this contradicts (\ref{eq3.1}).  Thus, we have $U_1'=U_1$.

Since $U_1'=U_1$ implies $U_2'=U_2$, and $U_2'\not=\emptyset$ in this case,
$G'$ cannot be an extremal graph, and so does $G$ since $\mathcal{N}(K_3, G')=\mathcal{N}(K_3, G)$.
Therefore, it is sufficient to show $G$ is an extremal graph only if $U_2=\emptyset$ in the case when $N(u')=N(u)$ for any $u'\in U_1-\{u\}$.



\vskip 2mm
Let  $X\subseteq G_u$ satisfy (\ref{eq3.1}).  By Lemma \ref{lem3}, $\Delta(G_u)=k-1$. Moreover, $G_u$ has the following structural properties.
\begin{Claim}\label{Claim1} Let $v\in V(C_i)$, where $C_i$ is some component of $G_u-X$.

(1) If $d_{G_u}(v)=k-1$,  then $N_{U_2}(v)$ is an independent set and $[v,U_2]\subseteq \mathcal{L}$.



(2) If $d_{C_i}(v)=k-2$, then $G[N_{U_2}(v)]$ is a star or a triangle, together with some isolated vertices. Moreover, if $v_1\in N_{U_2}(v)$ and $vv_1\in \mathcal{H}\cup \mathcal{M}$, then $v_1$ is the center of the star with at least 3 vertices, or lies on the triangle.
\end{Claim}
\noindent{\bf Proof.}  (1) Since $d_{G_u}(v)=k-1$ and $E(G_u)\subseteq \mathcal{H}$ by Lemma \ref{lem3}, we have $|\mathcal{H}(v)|\geq k-1$. Let
$X'\subseteq V(G_v)$ satisfy (\ref{eq3.1}). By Lemma \ref{lemma1},  $X'=\mathcal{H}(v)\subseteq V(G_u)$ which in turn implies $G_v-X'$ has no edges by (\ref{eq3.1}), and $\mathcal{M}(v)=\emptyset$, and so $[v,U_2]\subseteq \mathcal{L}$.

(2) Since $d_{C_i}(v)=k-2$ and $E(G_u)\subseteq \mathcal{H}$, we have $|\mathcal{H}(v)|\geq k-2$. Let
$X'\subseteq V(G_v)$ satisfy (\ref{eq3.1}). By Lemma \ref{lemma1},  $\mathcal{H}(v)\subseteq X'$ and so $|X'\cap V(G_u)|\geq k-2$. Because $|X'|\leq k-1$,  we have  $|X'\cap U_2|\leq 1$, which implies  $G_v-\mathcal{H}(v)$ has no $2K_2$ by (\ref{eq3.1}), that is, $G[N_{U_2}(v)]$ is a star or a triangle, together with some isolated vertices. Moreover, if $vv_1\in \mathcal{H}\cup \mathcal{M}$, then we have $d_{G_u}(v)=k-2$ by Lemma \ref{lemma1}.
Since $|N(vv_1)|\ge k$ and  we have $|N(vv_1)\cap U_2|\geq 2$, that is, $v_1$ has at least 2 neighbors in
$G[N_{U_2}(v)]$. Hence, $v_1$ is the center of the star or lies on a triangle in $G[N_{U_2}(v)]$.
 $\hfill\square$

\vskip 2mm
Let  $C$  be the largest component of $G_u-X$. Then
 $\pi(C)=(k-1,...,k-1,k-2)$ by Lemma \ref{lem3}.
Let $z\in V(C)$ be the unique vertex with $d_C(z)=k-2$ and $N_{C}(z)=\{z_1,\ldots,z_{k-2}\}$. Since $\Delta(G_u)=k-1$, we have $|[V(C),X]|\leq 1$ and $[V(C),X]\subseteq [z,X]$. In addition, we have the following.

\begin{Claim}\label{Claim2} $|\left(\mathcal{H}(z)\cup \mathcal{M}(z)\right)\cap U_2|\leq 1$.
\end{Claim}

\noindent{\bf Proof.}  Assume $v_1,v_2\in N_{U_2}(z)$ with $zv_1,zv_2\notin \mathcal{L}$. By Claim \ref{Claim1}(2),  $v_1v_2v$ is a triangle  in $N_{U_2}(z)$ for some $v$.
 Since $|N(zv_i)|\ge k$ for $i=1,2$, $\{z_1,\ldots,z_{k-2}\}\subseteq N(v_1)\cap N(v_2)$ which contradicts  Claim \ref{Claim1}(1) since $v_1,v_2\in N_{U_2}(z_1)$ and $v_1v_2\in E(G)$.  $\hfill\square$
 \vskip 2mm
 If $|\left(\mathcal{H}(z)\cup \mathcal{M}(z)\right)\cap U_2|=1$, we  assume $zv_1\in \mathcal{H}\cup \mathcal{M}$ and $N(zv_1)\cap U_2=\{v_2,...,v_t\}$. By Claim \ref{Claim1}(2), $G[\{v_1,...,v_t\}]$ is a $K_3$ or a star with center $v_1$. If $N(zv_1)\cap V(C)\not=\emptyset$,
let $N(zv_1)\cap V(C)=\{z_1,\ldots,z_{t'}\}$,  where $t'\leq k-2$.
\vskip 2mm
Let us consider the total weight in $\sum_{v\in V(C)}f(v)$ coming from the triangles not contained in $C$.
Since $[U_1,V(C)]\subseteq \mathcal{L}$, $U_1$ is an independent set, $|[V(C),X]|\leq 1$ and  $[X,U_2]\subseteq \mathcal{L}$ by Lemma \ref{lem3} and Claim \ref{Claim1}, the weight is contributed  by the triangles intersecting only with $U_2$. By Claims \ref{Claim1} and \ref{Claim2}, only $z$ is contained in some triangles with two vertices in $U_2$, and so the weight coming from such triangles is $\sum_{i\leq t}w(zv_1v_i,z)+\lambda(z)$, where $\lambda(z)=w(zv_2v_3,z)$ if $v_1v_2v_3$ is a triangle and $\lambda(z)=0$ otherwise. Furthermore, by Claims \ref{Claim1} and  \ref{Claim2}, for any triangle containing two vertices in $C$ and one vertex in $U_2$, only $w(zv_1z_i,z_i)\not= 0$ for $i\le t'$ in the case $zv_1\in \mathcal{H}\cup \mathcal{M}$. Therefore, the  weight is
 $$f_{U_2}(z)=\sum_{i\leq t}w(zv_1v_i,z)+\lambda(z)+\sum_{i=1}^{t'} w(zv_1z_i,z_i).$$

\begin{Claim}\label{Claim3} If $zv_1\notin \mathcal{H}$, then $f_{U_2}(z)\leq k-1$, and if $zv_1\in \mathcal{H}$, then $f_{U_2}(z)\leq k-1+\frac{t'}{2}\leq \frac{3k}{2}-2$ and $zv_1$ contributes at least $\frac{1}{2}(t'+1)$ to the loss of $v_1$.
\end{Claim}
\pf If $zv_1\in \mathcal{L}$, then $\sum_{i\leq t'}w(z_izv_1,z_i)=0$ and hence $f_{U_2}(z)=\sum_{i\leq t}w(zv_1v_i,z)+\lambda(z)\le max\{3,k-1\}$ by Lemma \ref{lemma2}.

If $zv_1\in \mathcal{M}$, then $|N(zv_1)|=(t-1)+t'\le 2k-2$. By the definition of $w(T,\cdot)$ and Claim \ref{Claim1}, we get
$f_{U_2}(z)\le \frac{1}{2}(t-1)+\lambda+\frac{t'}{2}.$ Note that $\lambda(z)=0$ if $(t-1)+t'=2k-2$, we have $f_{U_2}(z)\le k-1$.

 If $zv_1\in \mathcal{H}$, then $t>k$ and so $\lambda(z)=0$.  By Lemma \ref{lemma2}, $\sum_{i\leq t}w(zv_1v_i,z)\leq k-1$ and so $f_{U_2}(z)\leq k-1+\frac{t'}{2}\leq \frac{3k}{2}-2$. Moreover, since  the triangle $z_izv_1$ satisfies $zz_i\in \mathcal{H}$ for $i\le t'$,  by Observation \ref{O3}, the edge $zv_1$ contributes at least  $\frac{1}{2}(t'+1)$ to the loss of $v_1$.
 $\hfill\square$

\vskip 2mm
By Lemma \ref{lem3},  either $|C|=2k-1-2|X|\geq k+1$ or $|C|=k+1$. Moreover, since each edge of $G$ is covered by triangles, if $X=\emptyset$, then $G_u$ has no isolated vertices. That is,  $G_u=C$ or $C\cup K_{k-1}$ if $X=\emptyset$.
\vskip 2mm
We distinguish  the following two cases separately according to $|C|$.
\vskip 2mm
%
\vskip 2mm
\noindent{\bf Case 1. $|C|=2k-1-2|X|$}
\vskip 2mm
In this case, the structure of $G$ are shown in Figure 2, where the thick edges are in $\mathcal{H}$ and the thin edges are in $\mathcal{L}$.

\begin{picture}(70,70)(-52,130)
\centering
{\includegraphics[width=3.8in]{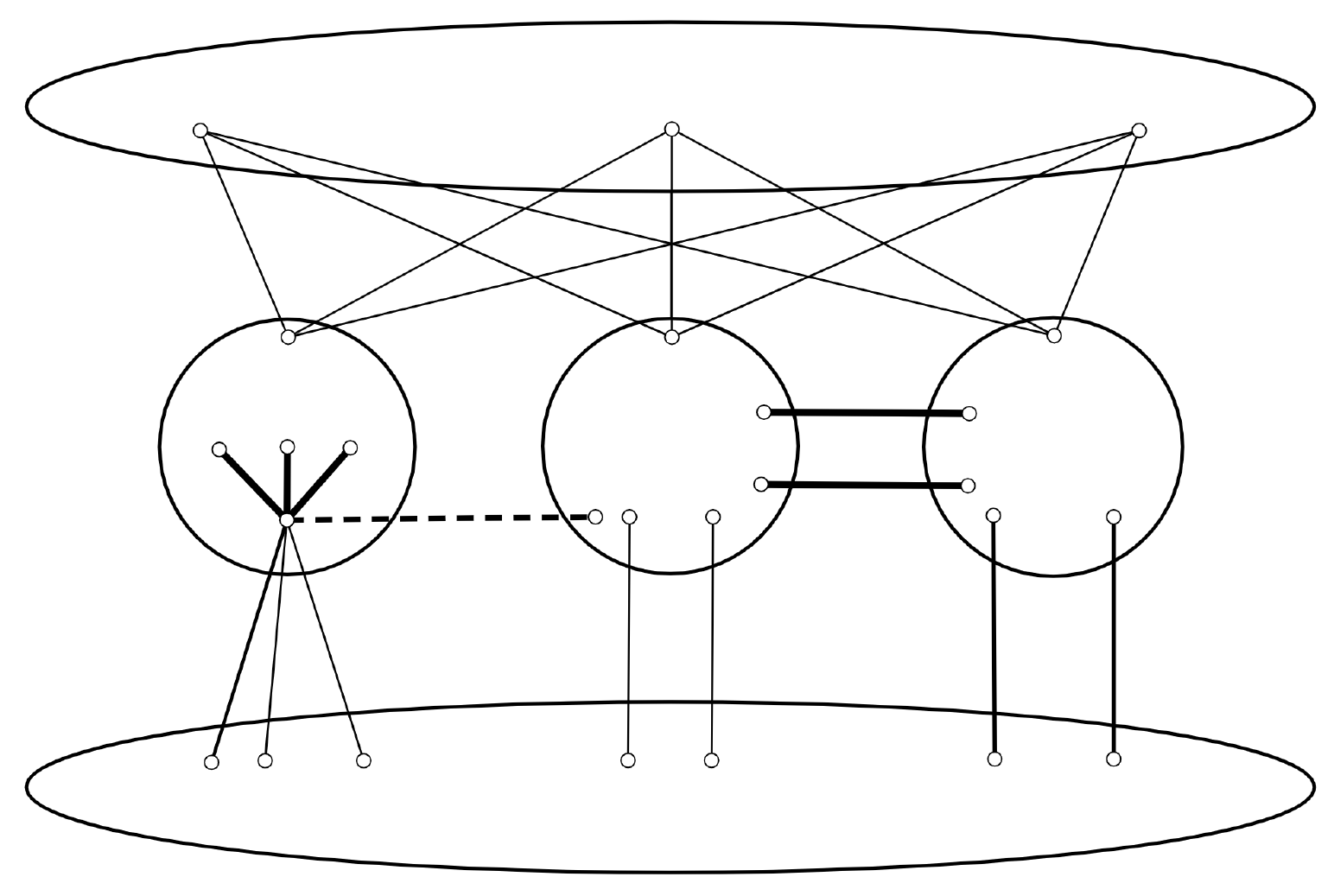}}
\put(-139,166){\makebox(3,3){$U_1$}}
\put(-216,105){\makebox(3,3){$C$}}
\put(-60,90){\makebox(3,3){$Y$}}
\put(-139,90){\makebox(3,3){$X$}}
\put(-223,73){\makebox(3,3){$z$}}
\put(-232,96){\makebox(3,3){$z_1$}}
\put(-204,96){\makebox(3,3){$z_{k-2}$}}
\put(-235,18){\makebox(3,3){$v_1$}}
\put(-220,18){\makebox(3,3){$v_2$}}
\put(-200,18){\makebox(3,3){$v_t$}}
\put(-139,14){\makebox(3,3){$U_2$}}

\put(-139,-20){\makebox(3,3){Figure 2. $N(u')=N(u)$ for any $u'\in U_1$}}
\end{picture}

\vskip 58mm

\vskip 2mm
\noindent{\bf Case 1.1} $X=\emptyset$.
\vskip 2mm
In this case,  $G_u=C$ and $|C|=2k-1$. If $U_2=\emptyset$, then since  $[U_1,G_u]\subseteq \mathcal{L}$ and $E(G_u)\subseteq \mathcal{H}$, by the definition of  $w(T,\cdot)$, we have
\[\begin{split}
\sum_{v\in V(G)}f(v)&=\sum_{v\in U_1}f(v)+\sum_{v\in V(C)}f(v)\\
&=\sum_{v\in U_1}f(v)+\mathcal{N}(K_3,C)=(n-2k+1)k\left(k-\frac{3}{2}\right)+\mathcal{N}(K_3,C).\\
\end{split}\]
Because $|C|=2k-1$ and $\pi(C)=(k-1,\ldots,k-1,k-2)$, by Theorem \ref{thm3},
$$\sum_{v\in V(G)}f(v)\le (n-2k+1)k\left(k-\frac{3}{2}\right)+2\binom{k-1}{3}+\binom{k/2}{2}+\binom{k/2-1}{2}$$
 and equality holds if and only if  $G=\bar K_{n-2k+1}+H_k$, and so the result follows. Therefore, we may assume that $U_2\not=\emptyset$.

 In this case, we will try to transfer the weight $f_{U_2}(z)$ to the vertices in $U_2$ such that $f(v)\le k(k-\frac{3}{2})$ still valid after transferring.

{\bf (1)} $zv_1\in \mathcal{H}$.

By Claim \ref{Claim3}, the edge $zv_1$ contributes at least $\frac{1}{2}(t'+1)$ to the loss of $v_1$.   Transfer the weight $\sum_{i\le t'}w(z_izv_1,z_i)$ to $v_1$ to cover the loss of $v_1$ caused by the edge $zv_1$ and the weight $\sum_{i\leq t}w(zv_1v_i,z)$ to $v_2,\ldots,v_t$($\lambda(z)=0$ in this case). After transferring, $f_{U_2}(z)=0$, $f(v)\le k\left(k-\frac{3}{2}\right)$ is still valid for any $v\in U_2$ and $f(v_1)\le k\left(k-\frac{3}{2}\right)-\frac{1}{2}$. Therefore,
\begin{equation}\label{eq5.1}
\sum_{v\in V(G)}f(v)< (n-2k+1)k\left(k-\frac{3}{2}\right)+\mathcal{N}(K_3,C).
\end{equation}

{\bf (2)} $zv_1\in\mathcal{M}$.

By the definition of $w(T,\cdot)$, we have $w(z_izv_1,z_i)=\frac{1}{2}$ for $i\leq t'$,
and $w(zv_1v_i,z)\leq \frac{1}{2}$ for $i\leq t$. Let $|U_2|-t=t''$.

 Suppose $t''\ge t'$.  Note that either $\lambda(z)=1$,  which implies $G[\{v_1,\ldots,v_{t}\}]$ is a triangle and $v_2v_3\in \mathcal{H}$, or $\lambda(z)\le \frac{1}{3}$. For the former case, $N(v_2v_3)-\{v_1,z\}\subset U_2-\{v_1,v_2,v_3\}$ and hence $|U_2-\{v_1,v_2,v_3\}|\geq 2k-3$ by Claim \ref{Claim1}(1), which means $t''>2t'$.
Thus  we can transfer the weight $w(z_izv_1,z_i)$ of $z_i$ ($i\leq t'$) to the vertices in $U_2-\{v_1,v_2,v_3\}$, the weight $\sum_{i\leq t}w(zv_1v_i,z)+\lambda(z)$ to the vertices $v_1,v_2,v_3$ and some others in $U_2-\{v_1,v_3,v_3\}$. For the latter case, we  transfer the weights $w(z_izv_1,z_i)$ to the vertices in $U_2-\{v_1,\ldots,v_t\}$ and  the weights $\sum_{i\leq t}w(zv_1v_i,z)+\lambda(z)$ to the vertices of $\{v_1,\ldots,v_t\}$. After the transferring, $f_{U_2}(z)=0$, $f(v)\le k\left(k-\frac{3}{2}\right)$ is still valid for any $v\in U_2$,  and (\ref{eq5.1}) still holds.

Assume $t''<t'\le k-2$. In this case, all edges of $G[\{v_1,\ldots,v_t\}]$ are in $\mathcal{L}$ for otherwise we have $t''\geq k-2$. Hence, $w(zv_1v_i,z)=\frac{1}{3}$ for $i\leq t$ and $\lambda(z)\le \frac{1}{3}$. Recalling $\{z_1,\ldots,z_{t'}\}\subseteq N(zv_1)$, by Claim \ref{Claim1}(1), we get $N(v_i)\cap \{z_1,\ldots,z_{t'}\}=\emptyset$ for $2\le i\le t$, and if $v_iz'v$ is a triangle such that $z'\in V(C)$ and $v\in U_2$, then $v_iz'v=v_izv_1$, or $G[\{v_1,...,v_t\}]=K_3$  and $v_iz'v=v_2zv_3$.
Thus, for $2\leq i\leq t$, we have
\[\begin{split}
f(v_i)&=\sum_{z'z''\subseteq C}w\left(v_iz'z'',v_i\right)+\sum_{v',v''\in U_2}w\left(v_iv'v'',v_i\right)+w\left(v_iv_1z,v_i\right)+\eta\\
 & \le\frac{1}{2}\Big(\left(2k-1-t'\right)(k-1)-t'(k-t')\Big)+\binom{t''+1}{2}+\frac{1}{3}+\frac{1}{3}\\
 &=k\left(k-\frac{3}{2}\right)-\frac{1}{2}\Big(t'(k+t'-1)-(t''+1)t''-1\Big)+\frac{2}{3},\\
\end{split}\]
where $\eta=w(v_2zv_3,v_2)$ or $w(v_2zv_3,v_3)$ if $v_2zv_3$ is a triangle, and $\eta=0$ otherwise. Because the total loss of the vertices in $U_2$ is at least
\[\begin{split}
&~\left(\frac{1}{2}\Big(t'(k+t'-1)-(t''+1)t''-1\Big)-\frac{2}{3}\right)(t-1)
+\frac{1}{2}+\frac{t''}{2},\\
\geq &~\sum_{i\leq t}w(zv_1v_i,z)+\lambda(z)+\sum_{i\le t'}w(z_izv_1,z_i)=\frac{t'}{2}+\frac{t}{3},\\
\end{split}\]
we can transfer these weights to vertices in $U_2$ and for $f_{U_2}(z)=0$, $f(v)\le k\left(k-\frac{3}{2}\right)$ is still valid for any $v\in U_2$,  and (\ref{eq5.1}) still holds.


{\bf (3)}  $zv_1\in \mathcal{L}$.

In this case, we have $w(z_izv_1,z_i)=0$ for $i\leq t'$ by the definition of $w(T,\cdot)$.
Since $zv_1\in \mathcal{L}$ and $\{v_2,...,v_t\}\subseteq N(zv_1)$, we have $t\leq k-1$. If $G[\{v_1,...,v_t\}]$ has an edge $v_iv_j \in \mathcal{H}$,
then since $N(v_iv_j)-\{z,v_1\}\subseteq U_2-\{v_1,...,v_t\}$ by Claim \ref{Claim1}(1), we have $|U_2-\{v_1,...,v_t\}|\geq 2k-3>k-1=t$, which implies $|U_2|> 2t$.  If $G[\{v_1,...,v_t\}]$ has no edge in $\mathcal{H}$, then $w(zv_1v_i,z)\le \frac{1}{3}$ and $\lambda(z)\le \frac{1}{3}$. Thus, by Lemma \ref{lem3}, we can we transfer the weight
$\sum_{i\leq t}w(zv_1v_i,z)+\lambda(z)$  to the vertices of $U_2$ in the former case and to the vertices of $\{v_1,...,v_t\}$ in the latter case. So,
$f_{U_2}(z)=0$, $f(v)\le k\left(k-\frac{3}{2}\right)$ is still valid for any $v\in U_2$,  and (\ref{eq5.1}) still holds.

Thus, Theorem \ref{thm3} and (\ref{eq5.1}) hold.

\vskip 3mm
\noindent{\bf Case 1.2} $X\neq \emptyset$.
\vskip 2mm
Let $X=\{x_1,\ldots,x_s\}~~\text{and}~~Y=V(G_u)-V(C)-X$. By Lemma \ref{lem3}, both $X$ and $Y$ are independent sets. Moreover, since $|C|=2k-1-2|X|=2k-1-2s\geq k+1$, we get that $s\leq \frac{k}{2}-1$.

By Claim \ref{Claim3},  $f_{U_2}(z)\le \frac{3k}{2}-2$.
Moreover, if $\big|[V(C),X]\big|=1$, then $d_{G_u}(z)=k-1$. By Claim \ref{Claim1}(1), $zv_1\in \mathcal{L}$ and $N_{U_2}(z)$ is an independent set. Therefore,
\begin{equation}\label{eq5-4}
f_{U_2}(z)=0~\text{if}~\big|[V(C),X]\big|=1.
\end{equation}

Consider the total loss of all the vertices in $Y$ contributed by the edges in $[X,Y]$.
Let $xy$ be any edge with $x\in X$ and $y\in Y$,
and $X_y\subseteq V(G_y)$ satisfy (\ref{eq3.1}).
Because $xy\in \mathcal{H}$, by Lemma \ref{lemma1}, $x\in \mathcal{H}(y)\subseteq X_y$. Since $X$ and $Y$ are independent sets,  $N(xy)\subseteq U_1\cup U_2$. Thus, if $xyv$ is a triangle, then  $xv\in \mathcal{L}$ by Claim \ref{Claim1}(1), which implies $w(xyv,y)=0$, and so $\sum_{v\in N(xy)}w(xyv,y)=0$. By Observation \ref{O2}, the edge $xy$ contributes $k-1$ to the loss of $y$.
Therefore, the total loss of all the vertices in $Y$, contributed by the edges in $[X,Y]$, is at least  $\big|[X,Y]\big|\cdot (k-1)$.

Since $X$ and $Y$ are independent sets, $N_{U_2}(x)$ is an independent set by Claim \ref{Claim1}(1),  $\big|[V(C),X]\big|\leq 1$ and $[U_1,G_u]\subseteq \mathcal{L}$. So,

$$f(x)=\sum_{y\in N_Y(x)}\sum_{v\in N(xy)\cap U_2}w(xyv,x)~~\text{for any}~ x\in X.$$
We try to transfer the weights $w(xyv,x)$ of $x\in X$ to the vertices $y,v$, such that the new weight of $x$ is 0, and that of
each other vertex remains no more than $k\big(k-\frac{3}{2}\big)$.

Fix an edge $yv$ and let $N(yv)\cap X=\{x_1,\ldots,x_{s'}\}$. Then  $x_iy\in \mathcal{H}$ and $x_iv\in \mathcal{L}$ for $1\le i\le s'$ by the arguments above. If $yv\in \mathcal{L}$, then $w(x_iyv,x_i)=0$, and so there is nothing to transfer. If $yv\notin \mathcal{L}$, then $w(x_iyv,x_i)=\frac{1}{2}$. Let $X' \subseteq V(G_v)$  satisfy (\ref{eq3.1}).

If $yv\in \mathcal{H}$, then since $yv,x_iy\in \mathcal{H}$, by Observation \ref{O3}, the edge $yv$  contributes at least $\frac{1}{2}(s'+1)$ to the loss of $v$, and so we can transfer the weight $\sum_{i=1}^{s'}w(x_iyv,x_i)=\frac{s'}{2}$ to $v$ to cover the loss caused by the edge $yv$.

Suppose $yv\in \mathcal{M}$. Note that $\sum_{v'\in N(x_iv)}w(x_ivv',v)\leq (k-1)-\frac{1}{2}$ by (\ref{eq3.4}) because $w(x_iyv,v)=\frac{1}{2}$.
If $x_i\in X'$, then by Observation \ref{O2}, the edge $x_iv$  contributes  $\frac{1}{2}$ to the loss of $v$. So we can transfer the weight $w(x_iyv,x_i)$ to $v$ to cover the loss contributed by the edge $vx_i$.
If $y\in X'$, then by (\ref{eq3.3}) and Observation \ref{O2}, the edge $yv$ contributes $\frac{k}{6}$ to the loss of $v$, and  $\frac{k}{12}$ to the loss of $y$
by Observation \ref{O4}, which means $yv$ contributes at least $\frac{k}{4}$ to the total loss of $y$ and $v$.
Recalling $s'\le s\le \frac{k}{2}-1$, we can transfer the weight $\sum _{i=1}^{s'}w(x_iyv,x_i)=\frac{s'}{2}$ to $y,v$ to cover the loss contributed by the edge $yv$. If neither $x_i\in X'$ nor $y\in X'$, then the edge $x_iy$ lies in  some component $C'$ of $G_v-X'$.
Remember $\sum_{v'\in N(x_iv)}w(x_ivv',v)\leq (k-1)-\frac{1}{2}$, the edge $x_iv$ contributes at least $\frac{1}{4}$ to the loss of $v$. Since $yv\in \mathcal{M}$, by Observation \ref{O4},  it contributes at least $\frac{k}{6}$ to the total loss of $y,v$. Thus, the total loss of $y$ and $v$ is at least $\frac{s'}{4}+\frac{k}{6}>\frac{s'}{2}$, and so we can transfer the weight $\sum _{i=1}^{s'}w(x_iyv,x_i)=\frac{s'}{2}$ to $y,v$ such that the weights of $y,v$ are still no more than $k(k-\frac{3}{2})$.

\vskip 2mm
If $zv_1\in \mathcal{H}$, then by Claim \ref{Claim3}, we can transfer $\frac{t'}{2}$ from $f_{U_2}(z)$ to $v_1$ to cover the loss of $v_1$ caused by $zv_1$ such that $f_{U_2}(z)\leq k-1$. After this possible transferring, we always have $f_{U_2}(z)\leq k-1$. Thus, recalling the total loss of all the vertices in $Y$ contributed by the edges in $[X,Y]$ is at least  $\big|[X,Y]\big|\cdot (k-1)$, (\ref{eq5-4}) and $f(x)=0$ for each $x\in X$ after transferring and Theorem  \ref{thm4}, the total weight of $G$ is

\[\begin{split}
\sum_{v\in V(G)}f(v)&\le (n-2k+1+s)k\left(k-\frac{3}{2}\right)+\mathcal{N}(K_3,C)+\left(k-1\right) -s(k-1)^2\\
&\le (n-2k+1)k\left(k-\frac{3}{2}\right)+2\binom{k-1}{3}+\left(\frac{k}{2}-1\right)^2\\
&~~~~+\frac{k^2}{4}+\left(-k^2+\frac{7k}{2}-\frac{11}{3}\right)s+(2k-2)s^2-\frac{4}{3}s^3~~  \\
&=(n-2k+1)k\left(k-\frac{3}{2}\right)+2\binom{k-1}{3}+\left(\frac{k}{2}-1\right)^2+\varphi(s,k).
\end{split}\]

Because $1\le s \le \frac{k}{2}-1$, after an easy calculation, we get $\varphi(k,s)<0$ except $\varphi(2,6)=1$, $\varphi(1,4)=3$ so the result holds  if $(s,k)\not=(2,6),(1,4)$.

Now, consider the two exceptions. Let $yv$ be any edge with $y\in Y$ and $v\in U_2$.

Suppose that $(s,k)=(2,6)$. Then $X=\{x_1,x_2\}$. Let $U=N_{U_2}(x_1)\cup N_{U_2}(x_2)$.
Assume  that $yv\in  \mathcal{H}\cup \mathcal{M}$. If both $x_1yv$ and $x_2yv$ are  triangles, then since $N_{U_2}(x_1)$ and $N_{U_2}(x_2)$ are independent sets by Claim \ref{Claim1},  $N(yv)-\{x_1,x_2\}\subseteq U_2-U$.  Clearly, $|U_2-U|\geq |N(yv)-\{x_1,x_2\}|\ge k-2$.
If $|N(yv)\cap \{x_1,x_2\}|\leq 1$ for any $yv\in \mathcal{H}\cup \mathcal{M}$, then since $k=6$,
$yv$ contributes at least $\frac{1}{2}$ to the loss of $y$ by Observations \ref{O3} and \ref{O4}, and so we can transfer $w(x_iyv,x_i)$ only to $y$ to cover the loss of $y$ contributed by $yv$ if $x_iyv$ is a triangle, such that  $f(x_i)=0$ for $i=1, 2$ after the possible transferring. Moreover, since $|N(yv)\cap U_2|\geq k-1$,  $|U_2|\geq k$. Thus, note that no weights are transferred to the vertices in $U_2-U$ in the former case and in $U_2$ in the latter case, $U_2$ has at least $k-2$ vertices whose weights are at most $k(k-\frac{3}{2})-\frac{1}{2}$ after transferring, which implies
\[\begin{split}
\sum_{v\in V(G)} f(v)&\le (n-2k+1)k\left(k-\frac{3}{2}\right)+2\binom{k-1}{3}+\left(\frac{k}{2}-1\right)^2+\varphi(s,k)-\frac{1}{2}(k-2)\\
&< (n-2k+1)k\left(k-\frac{3}{2}\right)+2\binom{k-1}{3}+\left(\frac{k}{2}-1\right)^2.\\
\end{split}\]
If $[Y,U_2] \subseteq \mathcal{L}$, then $f(x_i)=0$ for $i=1,2$ by Claim \ref{Claim1} and the definition of $w(T,\cdot)$.
Since $\varphi(2,6)=1$, if $|U_2|\geq 3$, then replace $\frac{1}{2}(k-2)$ with $\frac{1}{2}\cdot 3$ in the above inequality, we get the desired result. If $|U_2|\le 2$,  then $f(y)\leq 1$ for any $y\in Y$. It is easy to see the total weight of $G$ is less than the expected number.

Suppose that $(s,k)=(1,4)$. Then $X=\{x_1\}$.  We will transfer the weight $f(x_1)$ and  $f_{U_2}(z)$ to other vertices in a bit different way.  Note that $f(x_1)$ comes from the triangles $x_1yv$ with $yv\in \mathcal{H}\cup \mathcal{M}$.
Since $x_1y\in \mathcal{H}$ for any $y\in Y$, by Lemma \ref{lemma1}, $x_1\in X_y$.
If $yv\in \mathcal{H}$, then by Observation \ref{O3}, the edge $yv$ contributes $\frac{1}{2}$ to the loss of $y$. Assume that $yv\in \mathcal{M}$.
If $v\in X_{y}$, then by (\ref{eq3.3}) and Observation \ref{O2}, $yv$ contributes at least $\frac{k}{6}=\frac{2}{3}>\frac{1}{2}$  to the loss of $y$. If $v$ lies in a component $C'$ of $G_{y}-X_{y}$, then  $x_1\in X_y$ and (\ref{eq3.1}) imply $|C'|\leq 2k-3$, and so
$$\sum_{v'\in N(yv)\cap C'}w(yvv',y)\le \frac{1}{2}(k-2)+\frac{1}{3}\Big(2k-4-(k-2)\Big)=(k-1)-\left(\frac{k}{6}+\frac{2}{3}\right),$$
and so $yv$ contributes $\frac{1}{2}\left(\frac{k}{6}+\frac{2}{3}\right)>\frac{1}{2}$ to the loss of $y$ by Observation \ref{O4}. Therefore, the edge $yv$ contributes at least $\frac{1}{2}$ to the loss of $y$. Transfer the weight $w(x_1yv,x_1)=\frac{1}{2}$ to $y$ such that the new weight of $y$ is no more than $k\left(k-\frac{3}{2}\right)-(k-1)$, where the loss $k-1$ is contributed by the edge $x_1y$. Transfer the weight $f_{U_2}(z)$ to the vertices in $U_2$ in the same way used in Case 1.1. After the transferring, the weight of $G$ satisfies
$$\sum_{v\in V(G)}f(v)< (n-2k+1)k\left(k-\frac{3}{2}\right)+2\binom{k-1}{3}+\left(\frac{k}{2}-1\right)^2,$$
and so the proof of Case 1 is complete.

\vskip 3mm
\noindent{\bf Case 2.} $G_u=C_1\cup K_{k-1}$
\vskip 2mm
In this case, $X=\emptyset$ and $G_u=C\cup K_{k-1}$. Since $\pi(C)=(k-1,...,k-1,k-2)$ by Lemma \ref{lem3} and $|C|=k+1$,  $C$ is the complement of $\frac{1}{2}(k-2)K_2\cup P_3$, and so we have
$$\mathcal{N}(K_3,G_u)=\binom{k-1}{3}+\binom{k+1}{3}-\frac{1}{2}(k-2)(k-1)-2(k-2)-1.$$

Set $V(K_{k-1})=\{p_1,\ldots,p_{k-1}\}$. We first discuss some properties of these vertices.

\begin{Claim}\label{Claim4}
If $q_1,q_2\in U_2$ such that $p_iq_1,p_iq_2\in \mathcal{H}\cup\mathcal{M}$, then \\
(1) $G[N_{U_2}(p_i)]$ consists of a triangle $q_1q_2q_3$ and some isolated vertices;\\
(2) $|N(p_iq_1)|=|N(p_iq_2)|=k$  and $\{p_1,...,p_{k-1}\}\subseteq N(q_1)\cap N(q_2)$;\\
(3) If $v\in U_2$ such that $p_jv\in \mathcal{H}\cup \mathcal{M}$, then $v\in \{q_1,q_2,q_3\}$. Moreover, if $v=q_3$, then $|N(p_jq_s)|=k$
for $1\leq j\leq k-1$ and $1\leq s\leq 3$.
\end{Claim}
\pf
 (1) Since $d_{G_u}(p_i)=k-2$ and $p_iq_1,p_iq_2\in \mathcal{H}\cup\mathcal{M}$,  by Claim \ref{Claim1}(2), $G[N_{U_2}(p_i)]$ is  a triangle containing $q_1,q_2$, say $q_1q_2q_3$, together with some isolated vertices.

(2) Since $|N(p_iq_1)|\ge k$ and $N(p_iq_1)\cap U_2=\{q_2,q_3\}$ by (1), $\{p_1,...,p_{k-1}\}\subseteq N(q_1)$.  By the symmetry of $q_1$ and $q_2$, $\{p_1,...,p_{k-1}\}\subseteq N(q_2)$, and so the result follows.

(3) Suppose $v\notin \{q_1,q_2,q_3\}$. By Claim \ref{Claim1}(2) and (2), $G[N_{U_2}(p_j)]$ is a triangle $vq_1q_2$ together with some isolated vertices.
Because $|N(p_jv)|\ge k$ and $N(p_jv)\cap U_2=\{q_1,q_2\}$, we have $\{p_1,...,p_{k-1}\}\subseteq N(v)$. Thus, $v\in N(p_iq_1)$ and hence $|N(p_iq_1)|\geq k+1$ which contradicts (2). Therefore, $v\in \{q_1,q_2,q_3\}$.
Moreover, if $v=q_3$, then since $|N(p_jq_3)|\ge k$ and $N(p_jq_3)\cap U_2=\{q_1,q_2\}$, we can deduce $\{p_1,...,p_{k-1}\}\subseteq N(q_3)$.
 After an easy check, we get that $|N(p_jq_s)|=k$ for $1\leq j\leq k-1$ and $1\leq s\leq 3$.  $\hfill\square$

\vskip 2mm

By Lemma \ref{lemma1}, we have $|(\mathcal{H}(p_i)\cup \mathcal{M}(p_i))\cap U_2|\leq 3$ for all $1\leq i\leq k-1$. Suppose $|(\mathcal{H}(p_i)\cup \mathcal{M}(p_i))\cap U_2|=3$ for some $i$. By Claim \ref{Claim4}, $G[N_{U_2}(p_i)]$ is a triangle $q_1q_2q_3$ together with some isolated vertices and $p_iq_1,p_iq_2,p_iq_3\in \mathcal{M}$.  Furthermore,  we have $p_jq_s\in E(G)$ and $|N(p_jq_s)|=k$ for $1\leq j\leq k-1$ and $1\leq s\leq 3$, and $q_1q_2,q_1q_3,q_2q_3\in \mathcal{H}\cup \mathcal{M}$. Because
$$\sum_{v\in N(p_iq_s)} w(p_iq_sv,q_s)=\frac{k}{3}\le \frac{k}{2}=k-1-\left(\frac{k}{2}-1\right),$$
the edge $p_iq_s$ contributes at least $\frac{k}{4}-\frac{1}{2}$ to the loss of $q_s$ by Observation \ref{O4}.
%
Hence, all the edges $p_jq_s$, $1\leq j\leq k-1$ and $1\leq s\leq 3$, contribute  at least
$3(k-1)\left(\frac{k}{4}-\frac{1}{2}\right)$ to the total loss of $q_1$, $q_2$ and $q_3$. On the other hand, note that $G[N_{U_2}(p_j)]$ is the triangle $q_1q_2q_3$ together with some isolated vertices by Claim \ref{Claim4}(1) and $[p_j,U_2-\{q_1,q_2,q_3\}]\subseteq \mathcal{L}$ by Claim \ref{Claim4}(3). So, the weight of $p_j$ contributed by the triangles not in $K_{k-1}$ is
$$\sum_{1\leq r<s\leq 3} w(p_iq_rq_s,p_i)+\sum_{p_j\not=p_i}\sum_ {1\leq s\leq 3}w(p_ip_jq_s,p_i)=3\cdot\frac{1}{3}+3(k-2)\cdot\frac{1}{3}=k-1.$$
By Claim \ref{Claim3},
$f_{U_2}(z)\leq \frac{3k}{2}-2$. Therefore, the total weight of $G$ is at most
\[\begin{split}
& (n-2k)k\left(k-\frac{3}{2}\right)+\mathcal{N}(K_3,G_u)+\left(\frac{3k}{2}-2\right)+(k-1)^2-3(k-1)\left(\frac{k}{4}-\frac{1}{2}\right)\\
<&(n-2k+1)k\left(k-\frac{3}{2}\right)+2\binom{k-1}{3}+\left(\frac{k}{2}-1\right)^2,
\end{split}\]
a contradiction. So we assume that $|(\mathcal{H}(p_i)\cup \mathcal{M}(p_i))\cap U_2|\le 2$ for $1\leq i\leq k-1$.
\vskip 2mm

Fix $p_i$ and let $N_{U_2}(p_i)=\{q_1,\ldots, q_r\}$. By
Claim \ref{Claim1}, we may assume that $G[N_{U_2}(p_i)]$ is a star with the center $q_1$ or a triangle $q_1q_2q_3$, and some isolated vertices.
Let
$$f_{U_2}(p_i)=\sum_{v,v'\in U_2} w(p_ivv',p_i)+\sum_{p_j\not=p_i} \sum_{v\in N(p_ip_j)\cap U_2}w(p_ip_jv,p_j).$$
It is clear that $\sum_{i=1}^{k-1}f_{U_2}(p_i)$ is the total weight of $V(K_{k-1})$ contributed by the triangles not contained in $K_{k-1}$. We will complete the proof by showing that $f_{U_2}(p_i)\leq \frac{3k}{4}-\frac{1}{2}$ and $f_{U_2}(z)<\frac{3k}{4}-\frac{1}{2}$, after some appropriate weight transferring.
\vskip 2mm
 If  $|(\mathcal{H}(p_i)\cup \mathcal{M}(p_i))\cap U_2|=2$, say  $p_iq_1,p_iq_2\notin \mathcal{L}$, then by Claim \ref{Claim4}, $N_{U_2}(p_i)$ is a triangle $q_1q_2q_3$ together with some isolated vertices, $|N(p_iq_1)|=|N(p_iq_2)|=k$ and so
\[\begin{split}
&f_{U_2}(p_i)=\sum_{v,v'\in \{q_1,q_2,q_3\}}w(p_ivv',p_i)+\sum_{p_j\not=p_i}w(p_ip_jq_1,p_j)+\sum_{p_j\not=p_i}w(p_ip_jq_2,p_j).
\end{split}\]
Since
$$\sum_{v\in N(p_iq_s)} w(p_iq_sv,q_s)\le \frac{k}{2}=(k-1)-\left(\frac{k}{2}-1\right)~\text{for}~s=1,2,$$
the edge $p_iq_s$ contributes at least $\frac{k}{4}-\frac{1}{2}$ to the loss of $q_s$ for $s=1,2$ by Observation \ref{O4}.
Because $\sum_{p_j\not=p_i}w(p_ip_jq_s,p_j)\le \frac{1}{2}(k-2)$, transfer $\frac{k}{4}-\frac{1}{2}$ from $\sum_{p_j\not=p_i}w(p_ip_jq_s,p_j)$ to $q_s$ for $s=1,2$. After transferring, we have
\begin{equation}\label{eq43}
f_{U_2}(p_i) \leq \frac{3}{2}+2\cdot\frac{1}{2}(k-2)-2\left(\frac{k}{4}-\frac{1}{2}\right)=\frac{k+1}{2}\leq \frac{3k}{4}-\frac{1}{2}.
\end{equation}
%

Now, let $|(\mathcal{H}(p_i)\cup \mathcal{M}(p_i))\cap U_2|\le1$. Assume $\mathcal{H}(p_i)\cup \mathcal{M}(p_i)\subseteq \{q_1\}$  by Claim \ref{Claim1}(2).
Because $p_ip_j\in \mathcal{H}$ and $p_iq_s\in \mathcal{L}$,  $w(p_ip_jq_s,p_j)=0$ for $2\le s\le r$ and hence
$$f_{U_2}(p_i)=\sum_{j=2}^{r} w(p_iq_1q_j,p_i)+\lambda(p_i)+\sum_{p_j\not=p_i} w(p_ip_jq_1,p_j),$$
where $\lambda(p_i)=w(p_iq_2q_3,p_i)$ if $q_1q_2q_3$ is a triangle and $\lambda(p_i)=0$ otherwise.  Using the same proof as that of Claim \ref{Claim3}, we have

\begin{Claim}\label{Claim5} $f_{U_2}(p_i)\leq k-1$ if $p_iq_1\notin \mathcal{H}$, and $f_{U_2}(p_i)\leq k-1+\frac{\ell}{2}\leq \frac{3k}{2}-2$  and $p_iq_1$ contributes at least $\frac{1}{2}(\ell+1)$ to the loss of $q_1$ if $p_iq_1\in \mathcal{H}$,
where $|N(p_iq_1)\cap V(K_{k-1})|=\ell.$
\end{Claim}
In order to show  $f_{U_2}(p_i)\leq \frac{3k}{4}-\frac{1}{2}$ in this case and $f_{U_2}(z)<\frac{3k}{4}-\frac{1}{2}$, we need to consider the structure of $G[U_2]$.

If $vv'\in \mathcal{M}$ is an edge in $G[U_2]$, then by Observation \ref{O4}, $vv'$ contributes $\frac{k}{12}$ to the loss of $v$ and $v'$, respectively, that is, $vv'$ contributes $\frac{k}{6}$ to the total loss of vertices in $U_2$. On the other hand, by Claims \ref{Claim3} and  \ref{Claim5}, we can transfer some weight from $f_{U_2}(z)$ and $f_{U_2}(p_i)$ to $v_1$ and $q_1$, respectively,  such that $f_{U_2}(z)\leq k-1$ and $f_{U_2}(p_i)\leq k-1$, and $f(v_1)\leq k\left(k-\frac{3}{2}\right)$ and $f(q_1)\leq k\left(k-\frac{3}{2}\right)$ still hold. This together with (\ref{eq43})  implies that after transferring some weights to the vertices in $U_2$, the total weight in $\sum_{v\in V(G_u)}f(v)$ coming from the triangles not in $G_u$ is at most  $k(k-1)$. Therefore,  if $G[U_2]$ has $\frac{3k}{2}-2$ edges in $\mathcal{M}$, then we have
\[\begin{split}
\sum_{v\in V(G)} f(v)&\le (n-2k)k\left(k-\frac{3}{2}\right)+\mathcal{N}(K_3,G_u)+k(k-1)-\left(\frac{3k}{2}-2\right)\cdot \frac{k}{6}\\
&<(n-2k+1)k\left(k-\frac{3}{2}\right)+2\binom{k-1}{3}+\left(\frac{k}{2}-1\right)^2,
\end{split}\]
a contradiction. Hence, $G[U_2]$ contains at most $\frac{3k}{2}-3$ edges in $\mathcal{M}$. Moreover, we have

\begin{Claim}\label{Claim6}
Let $qq'\in \mathcal{H}$ be an edge in $G[U_2]$. If $\mathcal{M}(q)\cap \mathcal{M}(q')\cap \{p_1,...,p_{k-1}\}=\emptyset$, then $qq'$ contributes $\frac{3k}{8}$ to the total loss of $q$ and $q'$. Furthermore, $G[U_2]$ contains at most $\frac{3k}{2}-2$ edges in $\mathcal{H}\cup \mathcal{M}$.
\end{Claim}
\pf By Claim \ref{Claim2} and the assumption, $\mathcal{M}(q)\cap \mathcal{M}(q')\cap \{z,p_1,...,p_{k-1}\}=\emptyset$. Noting that
$G[U_2]$ has at most $\frac{3k}{2}-3$ edges in $\mathcal{M}$, we have $|\mathcal{M}(q)|+|\mathcal{M}(q')|\le k+\frac{3k}{2}-3$. Thus,  by (\ref{eq3.2}) and Lemma \ref{lemma1}(ii), we get
\[\begin{split}
&~\sum_{v\in N(qq')} w(qq'v,q)+\sum_{v\in N(qq')} w(q'qv,q')\\
\le&~\frac{1}{2}\left(\big|\mathcal{H}(q)-\{q'\}\big|+\big|\mathcal{M}(q)\big|+               \big|\mathcal{H}(q')-\{q\}\big|+\big|\mathcal{M}(q')\big|\right)\le 2(k-1)-\frac{3k}{8},
\end{split}\]
and so the conclusion follows by Observation \ref{O2}.

 In addition, recall $|(\mathcal{H}(p_j)\cup \mathcal{M}(p_j))\cap U_2|\le 2$ for all $1\leq j\leq k-1$, by Claim \ref{Claim4}, $G[U_2]$ has at most one edge $q_1'q_2'$ such that
$\mathcal{M}(q_1')\cap \mathcal{M}(q_2')\cap \{p_1,...,p_{k-1}\}\not=\emptyset$.
For any other $\mathcal{H}$-edge $qq'$ in $G[U_2]$, $qq'$ contributes at least $\frac{3k}{8}>\frac{k}{6}$ to the total loss of $q$ and $q'$, which,  together with the possible edge $q_1'q_2'$,  implies $G[U_2]$ contains at most $\frac{3k}{2}-2$ edges in $\mathcal{H}\cup \mathcal{M}$.  $\hfill\square$
\vskip 2mm
Now, let us re-consider $f_{U_2}(z)$ and $f_i=f_{U_2}(p_i)$ based on Claim \ref{Claim6}. For convenience, let $a\in \{z,p_1,\ldots, p_{k-1}\}$, $N_{U_2}(a)=\{b_1,\ldots,b_m\}$ is a star with center $b_1$ or a triangle $b_1b_2b_3$ in $G[N_{U_2}(a)]$, $(\mathcal{H}(a)\cup \mathcal{M}(a))\cap U_2\subseteq \{b_1\}$ and $N(ab_1)\cap V(G_u)=\{a_1,\ldots,a_\ell\}$. Recall the expressions of  $f_{U_2}(z)$ and $f_i=f_{U_2}(p_i)$, we have
$$f_{U_2}(a)=\sum_{j=2}^{m} w(ab_1b_j,a)+\lambda(a)+\sum_{j=1}^{\ell} w(ab_1a_j,a_j).$$

If $ab_1\in\mathcal{H}$, then since $aa_i\in \mathcal{H}$, we can transfer the weight $w(ab_1a_i,a_i)$  to $b_1$ to cover the loss caused by the edge $ab_1$ by Observation \ref{O3}. For the weight $w(ab_1b_j,a)$, we have $w(ab_1b_j,a)\leq \frac{1}{2}$ with equality only if $b_1b_j\in \mathcal{H}\cup \mathcal{M}$. Thus, by Claim \ref{Claim6}, after transferring, we have
$$f_{U_2}(a)=\sum_{j=2}^{m}w(ab_1b_j,a)+\lambda(a)\leq max\left\{2,\frac{3k}{4}-1\right\}< \frac{3k}{4}-\frac{1}{2}.$$

Assume $ab_1\in \mathcal{M}$. Consider $w(ab_1b_j,a)$. If $b_1b_j\in \mathcal{H}$, then  $w(ab_1b_j,a)=\frac{1}{2}$, and $b_1b_j$ contributes $\frac{3k}{8}$ to the total loss of $b_1$ and $b_j$ by Claim \ref{Claim6}. Since $a\in \{z,p_1,\ldots, p_{k-1}\}$, there are at most $k$ such triangles, and so we can transfer $\frac{3}{8}$ of each $w(ab_1b_j,a)$ to $b_1$ and $b_j$ to cover the total loss of $b_1$ and $b_j$ caused by the edge $b_1b_j$. If $b_1b_j\in \mathcal{M}\cup \mathcal{L}$, then $w(ab_1b_j,a)=\frac{1}{3}$ and $ab_1$ contributes $\frac{k}{12}$ to the loss of $b_1$ by Observation \ref{O4}. Thus, we can transfer $\frac{1}{12}$ of each $w(ab_1b_j,a)$ to cover the loss of $b_1$ caused by the edge $ab_1$. After transferring, we have $w(ab_1b_j,a)\le \frac{1}{4}$ and so
$$f_{U_2}(a)\le \frac{1}{4}|N(ab_1)\cap U_2)|+\frac{1}{2}|N(ab_1)\cap V(G_u)|\le \frac{3k}{4}-1< \frac{3k}{4}-\frac{1}{2}.$$

If $ab_1\in \mathcal{L}$, then $w(ab_1a_j,a_j)=0$ for $1\leq i\leq \ell$.
If $b_1b_j\in \mathcal{H}$, then $w(ab_1b_j,a)=1$. By Claim \ref{Claim6},  we can transfer $\frac{3}{8}$ to cover the total loss of $b_1$ and $b_j$ caused by the edge $b_1b_j$, with at most one exceptional edge in $G[U_2]$. If $b_1b_j\in \mathcal{M}\cup \mathcal{L}$, then $w(ab_1b_j,a)=\frac{1}{3}$. After transferring, we have $w(ab_1b_j,a)\le \frac{5}{8}$ with at most one exception and so
$$f_{U_2}(a)=\sum_{j=2}^{m} w(ab_1b_j,a)+\lambda(a)+\le \frac{5}{8}(k-2)+1< \frac{3k}{4}-\frac{1}{2}.$$

By the three inequalities above, we have $f_{U_2}(z)<\frac{3k}{4}-\frac{1}{2}$. Moreover, combining the three inequalities with (\ref{eq43}), we have $f_{U_2}(p_i)\leq \frac{3k}{4}-\frac{1}{2}$. Thus, after appropriate weight transferring, we have $f_{U_2}(z)+\sum_{i=1}^{k-1}f_i<k(\frac{3k}{4}-\frac{1}{2})$. Hence, the total weight of $G$ is
\[\begin{split}
\sum_{v\in V(G)} f(v)&< (n-2k)k\left(k-\frac{3}{2}\right)+\mathcal{N}(K_3,G_u)+k\left(\frac{3k}{4}-\frac{1}{2}\right)\\
&\le (n-2k+1)k\left(k-\frac{3}{2}\right)+2\binom{k-1}{3}+\left(\frac{k}{2}-1\right)^2.
\end{split}\]
The proof of Theorem \ref{main} is complete.
$\hfill\blacksquare$

\section{Concluding Remarks}

Theorem \ref{main} is proved for $n\ge 4k^3$, but notice that the statement does not hold for small $n$. For example, take at most five disjoint copies of $K_{2k}$ then the number of copies $K_3$ is more than the extremal number in the theorem. It would be nice to determine the sharp bound for $n$ when this generalized Tur\'an number is correct.

It is natural to ask what happens if we count larger cliques. 
The third author \cite{gerbner3} showed that $ex(n,K_r,F_k)=O(n)$ for every $k$ and $r$, but the constants in the upper bound are large. We conjecture that the extremal graph for $ex(n,K_r,F_k)$ is still $\bar K_{n-v(H)}+H$, where $H$ is a graph with $V(H)=k-1, \Delta(H)=k-1$. 

Let $H^T$ be the graph obtained by replacing each edge of $H$ with a triangle, e.g., the friendship graph can be considered as a $S_k^T$. So it is also interesting to ask what if we replace each edge of any other graph $H$ with a triangle? For example, consider the extremal function $ex(n,K_3,P_k^T), ~ex(n,K_3,C_k^T)$.

We have determined the largest number of triangles in $F$-free graphs when $F$ is a friendship graph, but not when $F$ is an extended friendship graph. Alon and Shikhelman \cite{Alon} showed that in that case $c_1|V(F)|^2n\le \ex(n,K_3,F)\le c_2|V(F)|^2n$ for absolute constants $c_1$ and $c_2$. Better bounds were obtained for some forests, including exact results for stars \cite{chase}, paths \cite{Luo} and forests consisting only of path components of order different from 3 \cite{zc}.

\vskip 3mm
\noindent{\bf\large Acknowledgements}

This research was supported by NSFC under grant numbers  11871270, 12161141003 and 11931006 and by the National Research, Development and Innovation Office NKFIH, grants  K132696, SNN-135643, KH 130371, SNN 129364, FK 132060, and KKP-133819..


\end{document}